\documentclass[12pt]{article}

\usepackage{latexsym}
\newtheorem{theoreme}{{\bf Th\'eor\`eme}}[section]

\newtheorem{corollaire principal}[principal]{\bf Corollaire}
\newtheorem{proposition}[theoreme]{{\bf Proposition}}
\newtheorem{lemme}[theoreme]{{\bf Lemme}}
\newtheorem{corollaire}[theoreme]{{\bf Corollaire}}

\newenvironment{demonstration}{\noindent{\bf D\'emonstration
}}{\nolinebreak $\Box $\hspace{-2.15mm}\raisebox{1.25mm}{.} \medskip}

\newenvironment{demonstration du lemme}{\noindent{\bf D\'emonstration du lemme
}}{\nolinebreak $\Box $\hspace{-2.15mm}\raisebox{1.25mm}{.} \medskip}

\def\RR{{\bf R}}

\def\ZZ{{\bf Z}}

\begin{document}


\title{{\Large \bf  Dynamique du pseudo-groupe des isom\'etries locales sur une vari\'et\'e Lorentzienne
analytique de  dimension $3$}}

\author{{\normalsize
{\bf Sorin DUMITRESCU}}}
\date{Avril  2006}

\maketitle

\vspace{0.3cm}
{\normalsize{{\bf Dynamique du pseudo-groupe des isom\'etries locales sur une vari\'et\'e Lorentzienne analytique de  dimension $3$}}}

\
\vspace{0.1cm}

\noindent{\bf{ R\'esum\'e.}} 
Soit $(M,g)$ une vari\'et\'e lorentzienne analytique r\'eelle  de dimension $3$ compacte et connexe. Nous d\'emontrons que l'existence d'une orbite ouverte (non vide) du pseudo-groupe des isom\'etries locales  implique que la m\'etrique lorentzienne est localement homog\`ene (i.e. le pseudo-groupe des isom\'etries locales de $g$ agit transitivement sur $M$).

\section{Introduction}

La pr\'esence d'une structure g\'eom\'etrique sur une vari\'et\'e diff\'erentiable $M$ induit  une partition de la vari\'et\'e en classes d'\'equivalence : deux points
se trouvent dans la m\^eme classe d'\'equivalence s'ils sont reli\'es par un diff\'eomorphisme local qui pr\'eserve la structure g\'eom\'etrique.

A l'instar du cadre riemannien, il convient d'appeler   {\it isom\'etrie locale}  une telle application qui pr\'eserve la structure g\'eom\'etrique et orbites du pseudo-groupe
des isom\'etries locales les classes de la partition pr\'ec\'edente.

Dans~\cite{Gro} M. Gromov montre que cette partition est  tr\`es r\'eguli\`ere pour les {\it structures g\'eom\'etriques rigides} (par exemple, pour
les m\'etriques pseudo-riemanniennes ou les connexions affines). Plus pr\'ecis\'ement, il existe un ouvert dense de $M$ dans lequel les orbites du pseudo-groupe des isom\'etries 
locales sont des sous-vari\'et\'es ferm\'ees. En particulier, si une telle orbite est dense alors celle-ci est ouverte et, par cons\'equent, la structure g\'eom\'etrique en question
est {\it localement homog\`ene} sur un ouvert dense (i.e. le pseudo-groupe des isom\'etries locales agit transitivement sur un ouvert dense).

Le r\'esultat pr\'ec\'edent est connu dans la litt\'erature  sous le nom du th\'eor\`eme de l'orbite-ouverte et a \'et\'e comment\'e et appliqu\'e par des nombreux auteurs~\cite{BFZ}, \cite{Ben},\cite{BFL},\cite{D}, \cite{DG}.  Il a \'et\'e utilis\'e de mani\`ere essentielle pour la classification des structures qui m\'elangent la g\'eom\'etrie et la dynamique  comme
les flots d'Anosov de contact~\cite{BFL}, dans l'\'etude des actions de ''gros groupes'' (par exemple, les r\'eseaux des groupes semi-simples) ou encore dans
l'\'etude des vari\'et\'es lorentziennes compactes dont le groupe d'isom\'etrie est non compact. Dans tous ces cas, on montre  \`a un moment  de la preuve qu'une
certaine  structure g\'eom\'etrique rigide localement homog\`ene sur un ouvert dense (d'apr\`es le th\'eor\`eme de l'orbite-ouverte), l'est en fait sur toute la vari\'et\'e.

Le th\'eor\`eme principal de cet article doit \^etre vu comme un approfondissement dans ce sens  du th\'eor\`eme de l'orbite-ouverte de M. Gromov  dans le cas  particulier des vari\'et\'es  lorentziennes analytiques r\'eelles de dimension $3$.
 
\begin{theoreme}   \label{homogene} : Soit $(M,g)$ une vari\'et\'e  lorentzienne analytique r\'eelle de dimension $3$ compacte et connexe.  Si le pseudo-groupe des isom\'etries locales
de $g$ admet une orbite ouverte non vide dans $M$, alors celle-ci  est \'egale \`a $M$.
\end{theoreme}

Dans le cadre analytique le th\'eor\`eme de l'orbite-ouverte est pr\'ecis\'e dans~\cite{DG},\cite{Gro} sous la forme suivante : {\it en dehors d'un ensemble analytique compact (eventuellement
vide), les orbites du pseudo-groupe des isom\'etries locales sont les fibres d'une fibration analytique de rang constant.}  Avec ce th\'eor\`eme, notre hypoth\`ese d'existence
d'une orbite ouverte pour le pseudo-groupe des isom\'etries locales est donc \'equivalente avec l'existence d'une orbite ouverte et {\it dense.}

Nous pouvons  alors \'enoncer  le th\'eor\`eme~\ref{homogene} sous la forme \'equivalente suivante :

\begin{corollaire}

Si le pseudo-groupe des isom\'etries locales admet une orbite dont l'adh\'erence est d'int\'erieur non vide, alors celle-ci  est \'egale \`a $M$ (la m\'etrique $g$ est localement homog\`ene).
\end{corollaire}

Remarquons que nos hypoth\`eses impliquent  automatiquement que tout invariant scalaire de la m\'etrique lorentzienne (par exemple, les fonctions sym\'etriques
des courbures principales) est constant. En effet, un tel invariant doit \^etre une fonction analytique constante sur un ouvert de $M$ et donc partout. En particulier,
les $3$ courbures principales de $g$ (d\'efinies comme les valeurs propres de la courbure de Ricci par rapport \`a la m\'etrique lorentzienne) sont constantes sur $M$.

  Dans le cas de la dimension $2$ ceci serait suffisant pour conclure : dans ce cas la courbure sectionnelle est un invariant scalaire et le fait que cette courbure sectionnelle
  soit constante implique n\'ecessairement  que la m\'etrique est localement homog\`ene~\cite{Wo}. Ceci est \'egalement suffisant dans le cadre riemannien  o\`u 
  (gr\^ace \`a la compacit\'e du groupe orthogonal) les invariants scalaires suffisent pour s\'eparer les orbites du pseudo-groupe des isom\'etries locales (voir~\cite{Tri} pour une  version effective de cette propri\'et\'e). Le th\'eor\`eme  principal de cet article est donc \'egalement vrai (en toute dimension) dans le cadre riemannien. Nous donnons  une
  preuve rapide de ce r\'esultat bien connu  \`a la proposition~\ref{riemannien}.
  
  Il convient de remarquer qu'en dimension $3$  la courbure ne se r\'esume pas \`a un scalaire (c'est un tenseur) et la courbure sectionnelle est une fonction m\'eromorphe
(en g\'en\'eral non constante) d\'efinie sur la $2$-grasmannienne. Cette fonction admet en g\'en\'eral des p\^oles en dehors de l'ouvert form\'e par les  plans non-d\'eg\'en\'er\'es.
  
  Notre preuve utilise l'analyticit\'e de mani\`ere essentielle et notamment la propri\'et\'e suivante de prolongement d'isom\'etries locales : {\it tout point $m$ de $M$ poss\`ede
  un voisinage ouvert $U_{m}$ tel que toute isom\'etrie locale proche de l'identit\'e d\'efinie sur un ouvert connexe $U$ contenu dans $U_{m}$ se prolonge \`a $U_{m}$}. Ce
  ph\'enom\`ene a \'et\'e d\'ecouvert pour la premi\`ere fois par K. Nomizu~\cite{No}  dans le cadre des m\'etriques riemanniennes analytiques et \'etendu par la suite par A. Amores~\cite{Amo} et
  M. Gromov~\cite{Gro} aux structures rigides analytiques. 
  
  L'argument pr\'ec\'edent combin\'e avec le principe de monodromie permet de voir que sur les vari\'et\'es analytiques compactes et simplement connexes les isom\'etries locales proches
  de l'identit\'e se prolonge en des isom\'etries globales. C'est pr\'ecis\'ement cette technique  qui est utilis\'ee dans~\cite{D} pour montrer que le  groupe des isom\'etries d'une vari\'et\'e
  lorentzienne analytique compacte et simplement connexe est n\'ecessairement compact.
  
  Avec cette remarque le th\'eor\`eme~\ref{homogene} est naturellement compl\'et\'e par le 

\begin{corollaire} 
Si de plus $M$ est simplement connexe,  alors $(M,g)$ est la sph\`ere $S^3$ munie d'une m\'etrique lorentzienne invariante par  l'action par translations du groupe de
Lie $S^3$ sur lui-m\^eme.
\end{corollaire}

La justification du corollaire est la suivante : le pseudo-groupe des isom\'etries locales proches de l'identit\'e agit transitivement sur $M$~\cite{Ben} (voir la section suivante pour
plus de d\'etails) ce qui implique qu'il existe un groupe (de Lie) connexe $G$ d'isom\'etries (globalement d\'efinies) qui agit transitivement sur $M$. Il vient que $M$ s'identifie au quotient de ce groupe par le stabilisateur d'un point.  D'apr\`es le r\'esultat de~\cite{D}, $G$ est compact et donc $M$ est un quotient de deux groupes
de Lie compacts. Le stabilisateur d'un point s'identifie alors \`a un sous-groupe compact du groupe lin\'eaire $O(2,1)$ : il s'agit n\'ecessairement  de l'identit\'e ou d'un sous-groupe
\`a un param\`etre compact stabilisateur d'un vecteur de norme \'egale \`a $-1$ (voir la section suivante pour des pr\'ecisions).

Dans le premier cas  que  $M$ s'identifie avec l'unique groupe compact connexe et simplement connexe de dimension $3$, qui est la sph\`ere $S^3$ et la m\'etrique lorentzienne
$g$ est invariante par l'action de $S^3$ sur lui m\^eme. La m\'etrique lorentzienne s'obtient \`a partir d'une forme quadratique de signature $(2,1)$ sur l'alg\`ebre de Lie
de $S^3$ et qui est transport\'ee par les translations de $S^3$ sur lui-m\^eme.

Dans le deuxi\`eme cas $M$ est le quotient d'un groupe de Lie $G$ compact connexe de dimension $4$ par un sous-groupe isomorphe \`a $S^1$.  Il n'existe que deux
possibilit\'es pour $G$ : ou bien $S^3 \times S^1$, ou bien $S^1 \times S^1 \times S^1 \times S^1$. Le tore $S^1 \times S^1 \times S^1 \times S^1$ ne poss\`ede aucun
quotient simplement connexe non trivial. Il vient que $G$ est isomorphe \`a $S^3 \times S^1$ et  que $M$ est un quotient de $S^3 \times S^1$ par un sous-groupe \`a un param\`etre compact. Comme $M$ est suppos\'ee simplement  connexe, il vient que la projection du stabilisateur d'un point de $M$ dans $S^3 \times S^1$ est surjective sur le facteur $S^1$ et donc le premier facteur isomorphe \`a $S^3$ agit simplement transitivement sur $M$. Nous  sommes ramen\'es donc au cas pr\'ec\'edent.

\section{Orbites du pseudo-groupe des isom\'etries locales}

Dans toute la suite  $M$ d\'esigne une vari\'et\'e analytique  r\'eelle de dimension $3$ compacte et connexe, munie d'une m\'etrique lorentzienne $g$.

Rappelons qu'une {\it isom\'etrie locale} de $(M,g)$ est un diff\'eomorphisme local entre deux ouverts de $M$ qui pr\'eserve $g$. L'ensemble des isom\'etries locales forment
un pseudo-groupe pour la composition. Dans le cas o\`u le pseudo-groupe des isom\'etries locales agit transitivement sur $M$ (avec une unique orbite), on dit
que $g$ est {\it localement homog\`ene}; autrement dit, on a unicit\'e du mod\`ele local de $g$.

Un pas important dans l'\'etude des orbites du pseudo-groupe des isom\'etries locales a \'et\'e fait par  I. Singer qui a d\'emontr\'e dans~\cite{Si} le r\'esultat suivant. Nous \'enon\c cons le r\'esultat
pour les m\'etriques lorentziennes (car c'est l'objet de cet article), mais celui-ci a \'et\'e d\'emontr\'e par I. Singer pour les m\'etriques riemanniennes et g\'en\'eralis\'e seulement
par la suite par M. Gromov pour les structures rigides~\cite{DG},~\cite{Gro} : 

{\it Si $g$ est une m\'etrique lorentzienne sur une vari\'et\'e compacte $M$, alors il existe un entier $s$ tel que $g$ soit localement homog\`ene d\`es que le $s$-jet
de $g$ est le m\^eme en chaque point de $M$.}

Pour pr\'eciser cet \'enonc\'e construisons le $s$-jet de $g$ en suivant la m\'ethode adopt\'ee dans~\cite{DG}, \cite{Gro}. Pour cela remarquons que la pr\'esence
d'une m\'etrique lorentzienne $g$ sur une vari\'et\'e $M$ permet de r\'eduire le groupe structural du fibr\'e des rep\`eres (qui est un $GL(3, \RR)$-fibr\'e) au groupe
orthogonal $O(2,1)$ qui pr\'eserve la forme quadratique de signature $(2,1)$. Ce groupe poss\`ede quatre composantes connexes et pour simplicit\'e nous allons
travailler avec la composante connexe de l'identit\'e qui est $PSL(2, \RR)$ : il est possible de consid\'erer que le groupe structural du fibr\'e des rep\`eres $g$-orthonorm\'es
$R(M)$ est $PSL(2, \RR)$, quitte \`a consid\'erer un rev\^etement fini non ramifi\'e de $M$ (cette op\'eration laisse invariant l'\'enonc\'e du th\'eor\`eme~\ref{homogene}).

Pour se convaincre que $PSL(2, \RR)$ est  bien le  groupe orthogonal en question, faisons agir $SL(2, \RR)$ par changement de variables sur l'espace vectoriel des formes binaires 
homog\`enes de degr\'e $2$. Il s'agit de l'espace vectoriel des \'el\'ements de la forme $ax^2 + 2bxy +cy^2$, avec $a,b$ et $c$ r\'eels; l'action de $SL(2, \RR)$ transite par $PSL(2, \RR)$
et pr\'eserve le discriminant $b^2-ac$ qui est bien une forme quadratique de signature $(2,1)$. Il sera utile par la suite de se rappeler que par rapport \`a cette forme quadratique
le vecteur $2xy$ est unitaire,
que son plan orthogonal est engendr\'e par les vecteurs isotropes $x^2$ et $-y^2$ dont le produit vaut $1$. 

Rappelons \'egalement que les sous-groupes \`a un param\`etre de $PSL(2, \RR)$ sont conjugu\'es \`a :

\begin{enumerate}

\item 

un sous-groupe {\it elliptique} de la forme $\left(  \begin{array}{cc}
                                                                 cos T   &   sinT \\
                                                                 -sin T     &  cos T \\
                                                                 \end{array} \right)$ qui fixe un vecteur
de norme $-1$;

\item   un sous-groupe {\it unipotent}  $\left(  \begin{array}{cc}
                                                                 1   &   T \\
                                                                 0     &  1\\
                                                                 \end{array} \right)$  qui fixe  un vecteur isotrope;

\item   un sous-groupe {\it semi-simple}  $\left(  \begin{array}{cc}
                                                                 T   &   0\\
                                                                 0     &  T^{-1} \\
                                                                 \end{array} \right)$ qui fixe un vecteur de norme $1$.

\end{enumerate}

Retournons maintenant au concept de jet d'une m\'etrique lorentzienne. Pour consid\'erer des jets d'ordre fini il est n\'ecessaire d'avoir des syst\`emes de coordonn\'ees.
En pr\'esence d'une m\'etrique lorentzienne $g$ il convient de privil\'egier les coordonn\'ees exponentielles (dans lesquelles le $1$-jet de $g$ est celui de la m\'etrique standard 
$db^2-da \cdot dc$). Si on fixe un point $m$ de $M$, un syst\`eme de coordonn\'ees exponentielles centr\'e en $m$ est enti\`erement d\'etermin\'e par le choix d'une base $g$-orthonorm\'ee de l'espace tangent $T_{m}M$. L'ensemble des tous les syst\`emes de coordonn\'ees exponentielles sur $M$ n'est donc rien d'autre que le fibr\'e des rep\`eres
orthonorm\'es $R(M)$. Si on associe \`a chaque syst\`eme de coordonn\'ees exponentielles centr\'e en $m$, le $s$-jet de la m\'etrique $g$ en $m$, on construit
une application $g^{(s)}$ d\'efinie sur $R(M)$ et \`a valeurs dans l'espace affine des $s$-jets de m\'etriques lorentziennes \`a trois variables dont le $1$-jet vaut $db^2-da \cdot dc$.
Quand on change la base de l'espace tangent en $m$ (et donc le syst\`eme de coordonn\'ees exponentielles centr\'e en $m$), le $s$-jet de $g$ en $m$ dans le nouveau syst\`eme de coordonn\'ees change par une action lin\'eaire et alg\'ebrique de $PSL(2, \RR)$.
Comme l'action de $PSL(2, \RR)$ fixe le $1$-jet $db^2-da \cdot dc$ vu comme l'origine de l'espace affine, nous pouvons oublier cet origine et voir l'action de $PSL(2, \RR)$
comme \'etant une action lin\'eaire sur un espace vectoriel $V^{(s)}$, dont la dimension augmente vertigineusement quand $s$ augmente...

L'application $s$-jet de $g$ est donc une application analytique $PSL(2, \RR)$-\'equivariante d\'efinie sur le fibr\'e $R(M)$ et \`a valeurs dans l'espace vectoriel $V^{(s)}$
muni d'une action lin\'eaire et alg\'ebrique de $PSL(2, \RR)$. Ce type d'application s'interpr\`ete \'egalement comme une section du fibr\'e vectoriel obtenu \`a partir
du fibr\'e $R(M)$ via la repr\'esentation du groupe structural $PSL(2, \RR)$ sur $V^{(s)}.$

Le th\'eor\`eme de I. Singer affirme que $g$ est localement homog\`ene si pour un entier $s$ suffisamment grand l'image de $g^{(s)}$ est constitu\'ee d'une unique
orbite de $V^{(s)}$ sous l'action de $PSL(2, \RR)$. Plus pr\'ecis\'ement, il est montr\'e dans~\cite{DG}, \cite{Gro} que tout $1$-jet d'application de $M$ dans $M$
s'int\`egre en une isom\'etrie locale s'il pr\'eserve le $s$-jet de $g$ : deux points $m$ et $m'$ sont reli\'es par une isom\'etrie locale si le $s$-jet de $g$ est le m\^eme
en $m$ et en $m'$ (en tant que orbite de $V^{(s)}$ sous l'action de $PSL(2, \RR)$).

En particulier, nous venons de voir qu'une isom\'etrie locale de $g$ est enti\`erement d\'etermin\'ee (sur un ouvert connexe) par son $1$-jet en un point. C'est pr\'ecis\'ement
cette propri\'et\'e qui exprime que $g$ est une {\it structure rigide} (d'ordre 1)  \`a la Gromov. Ceci implique que le pseudo-groupe des isom\'etries locales est un pseudo-groupe
de Lie (de dimension finie) et que toute isom\'etrie proche de l'identit\'e s'obtient en int\`egrant un champ de vecteurs analytique dit {\it champ de Killing}.

Il sera utile par la suite de se rappeler que dans la cadre analytique tout point $m$ de $M$ admet un voisinage ouvert $U_{m}$ dans $M$ tel que tout champ de Killing
d\'efini dans un ouvert connexe $U$ contenu dans $U_{m}$ se prolonge \`a $U_{m}$. Cette propri\'et\'e implique que {\it la fibre du faisceau des germes de champs de Killing
de $g$ est une alg\`ebre de Lie de dimension finie qui ne d\'epend pas du point} \cite{Gro}. Sous les hypoth\`eses du th\'eor\`eme~\ref{homogene} cette alg\`ebre de Lie $\cal G$
doit agir transitivement sur un ouvert de $M$ et elle est donc de dimension au moins $3$ : le morphisme (d'espaces vectoriels) de sp\'ecialisation $\cal G$$ \to T_{m}M$ qui associe \`a un germe de champ de Killing $Z$ d\'efini au voisinage  du  point $m$ sa valeur $Z(m)$ au point $m$ doit \^etre surjectif aux points de l'ouvert en question. Notons
\'egalement que l'action de $\cal G$ \'etant libre sur le fibr\'e des rep\`eres $R(M)$, sa dimension est n\'ecessairement  inf\'erieure ou \'egale \`a $6$ (la dimension 
de $R(M)$).

Remarquons que le stabilisateur dans $PSL(2, \RR)$ du $s$-jet de la m\'etrique en un point $m$ de $M$ est un sous-groupe alg\'ebrique ferm\'e de $PSL(2, \RR)$
et que chaque sous-groupe \`a un param\`etre du stabilisateur fournit un sous-groupe \`a un param\`etre d'isom\'etries locales de $M$ qui fixent $m$ (et qui
sont lin\'earis\'ees en coordonn\'ees exponentielles). Autrement dit,  un stabilisateur de dimension strictement positive du $s$-jet de la m\'etrique en $m$ donne
un champ de Killing local au voisinage de $m$ et qui s'annule au point $m$. Inversement, en consid\'erant la diff\'erentielle du flot d'un champ de Killing
s'annulant en $m$ on construit un sous-groupe \`a un param\`etre du stabilisateur du $s$-jet de $g$ au point $m$ (a priori, nous n'avons pas 
tout le sous-groupe \`a un param\`etre, mais seulement un voisinage de l'identit\'e; la conclusion repose sur le fait que le stabilisateur doit \^etre un sous-groupe
alg\'ebrique de $PSL(2, \RR)$ et contient donc l'adh\'erence de Zariski de ce voisinage de l'identit\'e).

\subsection{Th\'eorie des invariants}

Nous commen\c cons par expliquer bri\`evement les id\'ees de la preuve  du th\'eor\`eme de M. Gromov sur la structure de la partition de $M$ en orbites du pseudo-groupe des
isom\'etries locales~\cite{Ben}, \cite{DG}, \cite{Gro}. L'argument clef de la preuve consiste en l'utilisation dans ce cadre du th\'eor\`eme de stratification de M. Rosenlicht~\cite{Ro2}
pour les actions alg\'ebriques.

Pour expliquer cela consid\'erons l'adh\'erence de Zariski $W$ de l'image de $g^{(s)}$ dans $V^{(s)}.$ L'ensemble $W$ est une vari\'et\'e alg\'ebrique affine $PSL(2, \RR)$-invariante.
Formellement nous pouvons descendre l'application $g^{(s)}$ sur $M$ et la voir comme une application d\'efinie sur $M$ (et non sur $R(M)$) \`a valeurs dans l'espace
des orbites $W //PSL(2, \RR)$. D'apr\`es le th\'eor\`eme de I. Singer,  les orbites du pseudo-groupe des isom\'etries locales de $g$ sont les fibres de cette application. En g\'en\'eral le quotient $W //PSL(2, \RR)$ n'est pas s\'epar\'e, mais le th\'eor\`eme de Rosenlicht~\cite{Ro2},\cite{Zi}  affirme que quitte \`a enlever un ferm\'e de Zariski $W'$ de $W$
le quotient $W \setminus W' / PSL(2, \RR)$ devient une vari\'et\'e  alg\'ebrique (s\'epar\'ee). Sur l'ouvert dense et invariant de $M$ (compl\'ementaire  d'un sous-ensemble analytique
eventuellement vide) qui est la pr\'eimage par $g^{(s)}$ de $W \setminus W'$  et o\`u le rang de la diff\'erentielle de $g^{(s)}$ est maximal les orbites du pseudo-groupe des isom\'etries 
locales sont donc les fibres d'une fibration de rang constant. 

Les hypoth\`eses du th\'eor\`eme~\ref{homogene} impliquent  que cette fibration est triviale (\`a valeurs dans un point) et que $g$ est localement homog\`ene   sur un ouvert dense.
D\'esignons par $M \setminus S$ l'ouvert dense maximal de $M$ sur lequel $g$ est localement homog\`ene, $S$ \'etant un sous-ensemble analytique compact de $M$
(pas n\'ecessairement  connexe). 

{\it Notre but est de montrer que $S$ et vide.}

L'image de la restriction de $g^{(s)}$ au fibr\'e des rep\`eres sur $M \setminus S$ est exactement une orbite $O$ de $V^{(s)}$. L'image de $g^{(s)}$ est donc enti\`erement contenue
dans l'adh\'erence $\bar O$ de $O$.  Il est classiquement connu que (pour les orbites des actions alg\'ebriques r\'eelles) $\bar O \setminus O$ est un ensemble semi-alg\'ebrique de dimension
strictement inf\'erieure  \`a la dimension de $O$ (contenu dans un ensemble alg\'ebrique de dimension strictement inf\'erieure \`a la dimension de $O$) et
donc  le ferm\'e invariant $\bar O \setminus O$ est constitu\'e d'orbites de dimension strictement inf\'erieure \`a la dimension de $O$~\cite{Hum}, \cite{Zi}. L'image de la restriction de $g^{(s)}$ au fibr\'e des rep\`eres au-dessus des
points de l'ensemble $S$ tombe n\'ecessairement  dans $\bar O \setminus O$.

Le lemme suivant sera utile pour la compr\'ehension de l'image de $g^{(s)}.$

\begin{lemme}  \label{orbites}

 L'image de $g^{(s)}$  ne contient aucune orbite de dimension $1$.
\end{lemme}

\begin{demonstration}

Consid\'erons le $s$-jet de la m\'etrique lorentzienne $g$ comme \'etant une application $g^{(s)} : R(M) \to V^{(s)}$ du fibr\'e des rep\`eres dans l'espace
vectoriel des $s$-jets de m\'etriques lorentziennes sur $\RR^3$ dont le $1$-jet  \`a l'origine est  $db^2-da \cdot dc$. Cette application est $PSL(2, \RR)$-\'equivariante et
la repr\'esentation $V^{(s)}$ de $PSL(2, \RR)$ se d\'ecompose en une somme directe de repr\'esentations irr\'eductibles. 

Rappelons la liste des
repr\'esentations irr\'eductibles de $PSL(2, \RR)$ :  pour chaque dimension impaire 
$2d+1$ sup\'erieure ou \'egale \`a trois il existe une unique repr\'esentation irr\'eductible $H_{2d+1}$ de $PSL(2, \RR)$ dans un espace vectoriel de dimension $2d+1$ qui peut 
\^etre vue comme la repr\'esentation induite sur les polyn\^omes homog\`enes \`a deux variables de degr\'e $2d$ \`a partir de la repr\'esentation canonique
sur $\RR^2$. Il est important de remarquer que la repr\'esentation {\it $H_{2d+1}$  ne contient aucune orbite de dimension $1$.} Plus pr\'ecis\'ement, toutes les orbites
sont de dimension $3$ except\'e  l'origine (qui est un point fixe) et les orbites des \'el\'ements qui sont  la  puissance $d$-\`eme d'une forme quadratique et qui sont de dimension $2$.
Si la forme quadratique est le carr\'e d'une forme lin\'eaire, il s'agit de l'orbite de l'\'el\'ement  $x^{(2d)}$  et le stabilisateur dans  $PSL(2, \RR)$ est  une extension
par un groupe fini du sous-groupe unipotent 
des matrices de la forme $\left(  \begin{array}{cc}
                                                                 1   &   T \\
                                                                 0     &  1 \\
                                                                 \end{array} \right)$. Si la forme quadratique est le produit de $2$-formes lin\'eaires distinctes sur $\RR^2$, la composante connexe de l'identit\'e
       du stabilisateur est un sous-groupe \`a un param\`etre semi-simple. Le dernier cas est fournit par une forme quadratique non d\'ecomposable sur $\RR$ et dans ce cas
       la composante connexe de l'identit\'e du stabilisateur est un sous-groupe \`a un param\`etre elliptique.

 Supposons maintenant que $O$ est une orbite de dimension strictement positive contenue dans l'image de $g^{(s)}$. Cette orbite admet alors une projection non triviale
 dans au moins une des repr\'esentations $H_{2d+1}$ qui entre dans la d\'ecomposition de $V^{(s)}$. La projection est donc de dimension sup\'erieure ou \'egale \`a deux,
 ce qui implique que la dimension de $O$ est sup\'erieure ou \'egale \`a $2$.
 \end{demonstration}
 
 Des r\'esultats simples de la th\'eorie des invariants permettent de conclure la preuve du th\'eor\`eme~\ref{homogene} sous l'hypoth\`ese simplificatrice suppl\'ementaire
 d'existence sur $M$ d'un champ de vecteurs analytique invariant par l'action du pseudo-groupe des isom\'etries locales de $M$. Nous pr\'esentons
 ce r\'esultat sous la forme d'un lemme qui sera utile plus loin.

 \begin{lemme}   \label{champ}
 
 S'il existe sur $M$ un champ de vecteurs analytique non singulier (ne s'annulant pas) $X$ de norme n\'egative ou nulle tel que le pseudo-groupe des isom\'etries locales de $g$ qui pr\'eservent $X$ admet une orbite ouverte  non vide dans $M$, alors celle-ci est \'egale \`a $M$.
 \end{lemme}
 
 Remarque : l'\'enonc\'e reste  \'egalement vrai pour un champ de vecteurs $X$ de norme positive, mais la preuve directe est plus difficile (voir la derni\`ere section de l'article). 
  
 \begin{demonstration}
 
Remarquons que la $g$-norme de $X$ est une fonction analytique sur $M$ constante sur l'ouvert sur lequelle le pseudo-groupe des isom\'etries locales agit transitivement et donc partout. Quitte \`a normaliser nous avons deux cas :
$g(X)$  vaut $0$ ou $-1$.

Supposons pour commencer que le champ de vecteurs $X$ est isotrope. Ceci revient \`a dire que le groupe structural du fibr\'e des rep\`eres $R(M)$ se r\'eduit
au groupe unipotent  des matrices de la forme $\left(  \begin{array}{cc}
                                                                 1   &   T \\
                                                                 0     &  1 \\
                                                                 \end{array} \right)$. Il faut penser qu'au lieu de consid\'erer  tous les rep\`eres $g$-orthonorm\'es,  on ne prend que ceux qui
                                                                 admettent $X$ comme vecteur base (ici $X$ joue le r\^ole de la forme quadratique $x^2$ dans la base $(x^2,2xy, -y^2)).$

       D\'esignons par $\cal R$ le fibr\'e principal de groupe structural unipotent et consid\'erons  le $s$-jet $g^{(s)}$ de $g$ comme \'etant  une application
                                                             $\left(  \begin{array}{cc}
                                                                 1   &   T \\
                                                                 0     &  1 \\
                                                                 \end{array} \right)$-\'equivariante de $\cal R$ dans l'espace vectoriel d'une repr\'esentation alg\'ebrique lin\'eaire du groupe unipotent.
                                                                 
        L'hypoth\`ese du lemme implique que l'image de $g^{(s)}$ est contenue dans  l'adh\'erence d'une orbite $O$. Or, par un th\'eor\`eme classique du \`a Konstant et Rosenlicht~\cite{Ro}  les orbites des action alg\'ebriques des groupes unipotents sont ferm\'ees dans la topologie de Zariski et
        donc aussi dans la topologie transcendentale (il importe de rappeler que le groupe additif unipotent $\RR$ est isomorphe 
        au groupe multiplicatif des r\'eels positifs en tant que groupe de Lie, mais non en tant que groupe alg\'ebrique). Ceci implique que l'image de $g^{(s)}$ est
        form\'ee d'une seule orbite sous l'action du groupe unipotent et donc que $g$ est localement homog\`ene. Remarquons que nous avons m\^eme montr\'e que tous
        les $1$-jets d'applications de $M$ dans $M$,  qui pr\'eservent $X$, pr\'eservent \'egalement le $s$-jet de la m\'etrique et s'int\`egrent donc en des isom\'etries locales. Le pseudo-groupe des
        isom\'etries locales pr\'eservant \`a la fois $g$ et $X$ agit donc transitivement sur $M$.
        
        Passons \`a pr\'esent au cas o\`u $X$ est de norme constante \'egale \`a $-1$. Dans ce cas le stabilisateur de $X$ est un sous-groupe \`a un param\`etre
        compact elliptique conjugu\'e dans $PSL(2, \RR)$ au groupe $\left(  \begin{array}{cc}
                                                                 cos T   &   sinT \\
                                                                 -sin T     &  cos T \\
                                                                 \end{array} \right)$.   Nous avons donc que le groupe structural du fibr\'e $R(M)$ se r\'eduit \`a un groupe compact. Comme les orbites
                                                                 d'un tel groupe (dans l'espace des $s$-jets de m\'etriques)  sont compactes et donc ferm\'ees les arguments du cas pr\'ec\'edent s'appliquent.
        \end{demonstration}
        
        Notons que  l'argument de la  preuve du lemme pr\'ec\'edent   s'applique \'egalement d\`es que le groupe structural du fibr\'e des rep\`eres se reduit \`a un groupe 
       compact. Dans le cas du  groupe orthogonal $O(n, \RR)$ nous avons alors la 
        
        \begin{proposition} \label{riemannien}
        
        Si $(M,g)$ est une vari\'et\'e riemannienne analytique connexe et compacte     telle que le pseudo-groupe d'isom\'etries locales de $M$ admet une orbite ouverte non vide,
        alors celle-ci est \'egale \`a $M$.
        \end{proposition}

\section{Dynamique des champs de Killing}
  
  R\'esumons bri\`evement la situation d\'ecrite dans la section pr\'ec\'edente. La m\'etrique lorentzienne $g$ est localement homog\`ene en dehors d'un ensemble
  analytique compact $S$ (admettant  eventuellement plusieures composantes connexes).   L'alg\`ebre de Lie $\cal G$ des germes de champs de Killing de $g$
  agit donc transitivement sur $M \setminus S$ : elle est de dimension au moins $3$.

  L'image $g^{(s)}(R( M \setminus S))$ du $s$-jet de $g$ sur l'ouvert $M \setminus S$
  est une  $PSL(2, \RR)$-orbite $O$ de l'espace vectoriel de jets $V^{(s)}$, tandis que l'image de $g^{(s)}$ du fibr\'e $R(M)$ restreint \`a $S$ est constitu\'ee d'un ensemble d'orbites contenues
  dans $\bar O \setminus O$ (car $R(M \setminus S)$ dense dans $R(M)$);  ces orbites sont de dimension strictement inf\'erieure \`a la dimension de $O.$ 
  
  {\it Notre but est de prouver que l'ensemble $S$ est n\'ecessairement vide. }
  
  Supposons par l'absurde que $S$ est non vide. La contradiction viendra d'une \'etude syst\`ematique  des
  diverses actions possibles de $\cal G$ au voisinage de l'ensemble $S$.
  
  Rappelons que l'action de $\cal G$ est libre sur le fibr\'e des rep\`eres orthonorm\'es $R(M)$ qui est de dimension $6$, ce qui implique que la dimension de $\cal G$ est
  au plus $6$. Nous d\'emontrons la

  \begin{proposition}  \label{dim3}
  L'alg\`ebre de Lie $\cal G$ est de dimension $3$.
  \end{proposition}
  
  \begin{demonstration}
  
  Remarquons d'abord que $\cal G$ ne peut pas \^etre de dimension $6$. En effet, la dimension $6$ correspond au cas de courbure sectionnelle constante (action transitive
  de $\cal G$ sur $R(M)$ et donc sur les $2$-plans non-d\'eg\'en\'er\'es) et dans ce cas
  $g$ est localement homog\`ene sur $M$ tout entier.
  
  \'Eliminons maintenant le cas o\`u $\cal G$ est de dimension $5$. Dans ce cas le groupe d'isotropie de $g$ sur l'ouvert $M \setminus S$ est de dimension $2$ car le morphisme
  d'\'evaluation $\cal G$$ \to T_{m}M$ en un point $m$ de $M \setminus S$ admet un noyau de dimension $2$. L'orbite $O$ poss\`ede alors un stabilisateur dans $PSL(2, \RR)$ de dimension $2$  et elle est, par cons\'equent, de dimension $1$.  Ce cas a \'et\'e exclus par le lemme~\ref{orbites}.

  Supposons \`a pr\'esent que $\cal G$ est de dimension $4$. Dans ce cas le groupe d'isotropie de $g$ sur $M \setminus S$ est de dimension $1$ et $O$ est de dimension $2$
  (poss\`ede un stabilisateur dans $PSL(2, \RR)$ de dimension $1$). Les orbites de $\bar O \setminus O$ sont alors de dimension strictement inf\'eriure \`a deux et d'apr\`es
  le lemme~\ref{orbites}  celles-ci sont de dimension $0$.  Nous avons donc que le groupe d'isotropie en un point de $S$ est de dimension $3$, isomorphe \`a $PSL(2, \RR)$.
  
   Fixons un point $s \in S$ et consid\'erons le morphisme d'\'evaluation qui \`a chaque germe de champs de Killing au voisinage de $s$
  associe sa valeur en $s$. On vient de constater que le noyau $\cal I$ de ce morphisme est de dimension $3$ ce qui implique que son image est un sous-espace vectoriel
  de dimension $1$ de $T_{s}M$ (qui est n\'ecessairement invariant  par l'action lin\'eaire du groupe d'isotropie $\cal I$).  La contradiction  vient du fait que l'action lin\'eaire    de $PSL(2, \RR)$ sur $T_{s}  M$ ne laisse stable  aucune  droite.
\end{demonstration}
 
 La proposition~\ref{dim3} implique le fait suivant :
 
 \begin{lemme}  \label{unimodulaire}
 
 i) L'alg\`ebre de Lie $\cal G$ est r\'esoluble et non-unimodulaire.
 
 ii) L'ouvert $M \setminus S$ admet une $(G,G)$-structure, o\`u $G$ est  l'unique groupe de Lie connexe et simplement connexe associ\'e \`a $\cal G$. En restriction \`a
 $M \setminus S$ la m\'etrique lorentzienne $g$ provient d'une m\'etrique lorentzienne sur $G$ invariante par translations.
  \end{lemme}
  
  Rappelons que l'existence d'une $(G,G)$-structure sur une vari\'et\'e \'equivaut \`a la donn\'ee d'un atlas \`a valeurs dans des ouverts de $G$ tel que les applications de changement de carte
  soient des restrictions de diff\'eomorphismes  de $G$ obtenus par  translations (\`a gauche). Tout objet g\'eom\'etrique d\'efini sur $G$ et invariant par les translations
  \`a gauche de $G$ sur lui-m\^eme fournit un objet g\'eom\'etrique de m\^eme nature sur l'ouvert $M \setminus S.$ Le lemme affirme que sur l'ouvert $M \setminus S$
  la m\'etrique lorentzienne est construite pr\'ecis\'ement de cette mani\`ere.  Il est utile d'observer que dans ce cas les champs de Killing locaux engendrent des translations (locales) \`a gauche dans $G$ et sont donc des champs de vecteurs invariants par les translations \`a droite dans $G$ : un tel champ de
  Killing n'est globalement d\'efini sur $M \setminus S$ que s'il est dans le centre de l'alg\`ebre de Lie $\cal G$.
  
  \begin{demonstration}
  
  La deuxi\`eme partie du lemme r\'esulte  du fait que le pseudo-groupe des isom\'etries locales de $g$ agit simplement transitivement sur $M \setminus S$ : sous l'effet 
  de cette action le voisinage de chaque point de $M \setminus S$ s'identifie \`a un ouvert de $G$ bien d\'efini \`a translation pr\`es.

  Passons maintenant \`a la preuve du premier point du lemme et supposons par l'absurde que $\cal G$ est unimodulaire~\cite{Mi}. Consid\'erons 
  $(K_{1}, K_{2}, K_{3})$ trois champs de Killing locaux
  lin\'eairement ind\'ependants : le caract\`ere unimodulaire du groupe implique alors que les translations \`a gauche pr\'eservent le volume de ces trois champs de vecteurs (le  volume \'etant calcul\'e par rapport \`a la m\'etrique lorentzienne $g$). Autrement dit, la fonction $vol(K_{1}, K_{2}, K_{3})$ est constante. 
  
  Fixons un point $s$ dans $S$ et consid\'erons un voisinage  $U$ de $s$ dans $M$ qui satisfait \`a la propri\'et\'e de prolongement de champs de Killing. Soient
  $(K_{1}, K_{2}, K_{3})$ trois champs de Killing lin\'eairement ind\'ependants dans un voisinage ouvert connexe contenu dans $U$ (ces trois champs existent
  car l'action de $\cal G$ est transitive sur $M \setminus S$). Ces trois champs se prolongent alors dans tout $U$ et  les vecteurs $(K_{1}(s), K_{2}(s), K_{3}(s))$
  ne sont pas libres, car inclus dans $T_{s}S$ \`a cause de l'invariance de $S$. La contradiction recherch\'ee vient du fait que la fonction analytique
  $vol(K_{1}, K_{2}, K_{3})$ est constante non nulle sur $U \setminus S$ et s'annule sur $S$.
  
  L'alg\`ebre de Lie $\cal G$ est n\'ecessairement r\'esoluble, car en dimension $3$ une alg\`ebre de Lie est ou bien unimodulaire, ou bien r\'esoluble~\cite{Mi}.
  \end{demonstration}
  
  \begin{proposition}  \label{surface}
  
  i) Les orbites de $\bar O \setminus O$ qui se trouvent dans l'image de $g^{(s)}$ sont toutes de dimension $2$.
  
  ii) Chaque composante connexe de $S$ est une surface lisse (sous-vari\'et\'e de codimension $1$ de $M$), d\'eg\'en\'er\'ee par rapport \`a $g$ et
  sur laquelle le pseudo-groupe des isom\'etries locales agit transitivement.
  
  iii) Chaque point de $S$   poss\`ede un voisinage ouvert $\Sigma$
  dans $S$ sur lequel il existe  un champ de vecteurs $X$ tangent \`a $\Sigma$ de $g$-norme constante \'egale \`a $0$  ou $1$ qui est pr\'es\'erv\'e par l'action de $\cal G$ (restreinte \`a $\Sigma$).
  \end{proposition}
  
  \begin{demonstration}
  
 i)  Rappelons  que les orbites de $\bar O \setminus O$ sont de dimension strictement inf\'erieure \`a la dimension de $O$ : celles-ci sont donc de
 dimension au plus deux. D'apr\`es le lemme~\ref{orbites} il ne peut y avoir dans l'image de $g^{(s)}$ des orbites de dimension $1$. Il reste \`a voir que dans l'image de $g^{(s)}$ il ne peut y avoir des orbites de dimension $0$. Pour cela notons que le stabilisateur d'une  orbite de dimension $0$ est isomorphe \`a $PSL(2, \RR)$
 et fournit donc une alg\`ebre d'isotropie de dimension $3$ isomorphe \`a l'alg\`ebre de Lie de $PSL(2, \RR)$. Nous avons donc que $\cal G$ est l'alg\`ebre de Lie
 de $PSL(2, \RR)$ qui est unimodulaire (car semi-simple) : contradiction avec le  lemme~\ref{unimodulaire}.
   
   ii) Remarquons que $\bar O \setminus O$ \'etant de dimension inf\'erieure ou \'egale \`a $2$, toute orbite de dimension $2$ contenue dans $\bar O \setminus O$ est un ouvert de
   $\bar O \setminus O$. Si $u$ est un point  de $S$ o\`u le $s$-jet de $g$ appartient \`a une orbite $O_{1}$ de dimension $2$ contenue dans $\bar O \setminus O$, alors
   par continuit\'e il existe un voisinage ouvert $\Sigma$ de $u$ dans $S$ o\`u le $s$-jet de $g$ appartient \`a $O_{1}$.  Ceci \'etant vrai au voisinage de chaque point $u$
   de $S$, il vient (par connexit\'e) que le $s$-jet de $g$ appartient \`a $O_{1}$ en chaque point de la composante connexe de $u$ dans $S$. Le pseudo-groupe
   des isom\'etries locales agit donc transitivement sur cette composante connexe.
   
   Une autre mani\`ere de le voir est de consid\'erer le morphisme d'\'evaluation au point $u$ de $S$. Comme le noyau est de dimension $1$ (car isotropie de dimension $1$),
   l'image est un sous-espace vectoriel de $T_{u}M$ de dimension $2$. L'ensemble $S$ \'etant invariant par l'action des champs de Killing nous avons que l'image du 
   morphisme d'\'evaluation en $u$ est incluse (et donc \'egale) \`a $T_{u}S$. Ceci implique que $S$ est de dimension $2$ et que le pseudo-groupe des isom\'etries locales agit
   transitivement sur un voisinage ouvert  de $u$ dans $S$. Par connexit\'e, l'orbite de $u$ sous l'action du pseudo-groupe des isom\'etries locales contient toute
   la composante connexe de $u$ dans $S$.
   
   Comme le pseudo-groupe des isom\'etries locales de $g$
   agit transitivement sur chaque composante connexe de $S$ tout en la  pr\'es\'ervant, il vient que chaque composante connexe de $S$ est une surface lisse ( sous-vari\'et\'e de codimension $1$ dans $M$).
   
   iii) L'\'enonc\'e qu'on veut d\'emontrer \'etant local, fixons  un point $u$ dans $S$ et supposons que $S$ est connexe. L'isotropie \'etant de dimension $1$ en les points de $S$, le stabilisateur $I_{1}$ de $O_{1}$ dans $PSL(2, \RR)$ est un sous-groupe alg\`ebrique de dimension $1$ : il s'agit d'une extension par un groupe fini d'un  sous-groupe \`a un param\`etre $I$ de $PSL(2, \RR)$. Le
    $s$-jet de $g$ au-dessus de $S$ est une application du fibr\'e des rep\`eres $R(M)$ restreint \`a $S$ dans l'orbite $PSL(2, \RR) /  I_{1}$. 
    
    En restriction \`a $S$ le groupe structural du fibr\'e des rep\`eres $R(M)$ se reduit au groupe $I_{1}$ et l'action de la restriction de $\cal G$ \`a $S$    
    pr\'eserve cette r\'eduction du groupe structural.
    
    Quitte \`a se restreindre  \`a un voisinage ouvert connexe  $\Sigma$ de $u$ dans $S$ suffisamment  petit, on peut consid\'erer que le groupe structural de $R(M)$ se r\'eduit \`a la composante connexe
    de l'identit\'e $I$ du groupe $I_{1}$. Comme l'espace homog\`ene  
    $PSL(2, \RR) /  I$
    repr\`esente l'ensemble des vecteurs non nuls de norme constante   de $\RR^{2,1}$ (\'egale \`a $0$, $1$ ou $-1$ selon que le sous-groupe \`a un param\`etre $I$
    est respectivement unipotent, semi-simple ou elliptique), 
   nous avons  en restriction \`a $\Sigma$ un champ de vecteurs $X$ de norme constante qui est
    pr\'es\'erv\'e par la restriction \`a $\Sigma$ de l'alg\`ebre de Lie des champs de Killing $\cal G$. On constate que pour tout $u \in \Sigma$ le vecteur $X(u)$ \'etant fix\'e en particulier par le champ de Killing  s'annulant en $u$ (qui provient du stabilisateur de $O_{1}$), toute la g\'eod\'esique issue de $u$ dans la direction de $X(u)$
    est form\'ee par des points fixes pour ce champ de Killing : cette g\'eod\'esique est n\'ecessairement  contenue dans $S$ (car l'isotropie est triviale sur $M \setminus S$). Nous avons en particulier que $X(u) \in T_{u} \Sigma$ et donc
    que $X$ est un champ de vecteurs tangent \`a $\Sigma$.
    
  Le champ $X$ ne peut pas \^etre de norme constante \'egale \`a $-1$ (ou de mani\`ere \'equivalente, le stabilisateur de $O_{1}$ ne peut pas \^etre elliptique). En effet,  dans 
  ce cas l'action (lin\'eaire) du groupe d'isotropie en un point $u$ de $\Sigma$ doit pr\'eserver  \`a la fois le plan $T_{u} \Sigma$ et le vecteur $X(u) \in T_{u} \Sigma$ de $g$-norme constante
  \'egale \`a $-1$. Ceci est impossible.
  
   Un raisonnement analogue montre que la m\'etrique lorentzienne $g$ est d\'eg\'en\'er\'ee sur $\Sigma$. En effet,  l'action du groupe d'isotropie en un point $u$ de $\Sigma$
   doit pr\'eserver le plan $T_{u}\Sigma$ et fixer le vecteur $X(u) \in T_{u} \Sigma$. Si le plan $T_{u} \Sigma$ n'est pas $g$-d\'eg\'en\'er\'e, l'action du groupe d'isotropie
   pr\'eserve n\'ecessairement un  vecteur $Z(u) \in T_{u}M$ de $g$-norme constante \'egale \`a $1$, orthogonal \`a $T_{u} \Sigma$. L'action de $I$ fixe donc tous les points de
   la g\'eod\'esique issue de $u$ dans la direction $Z(u)$. Comme ces points se trouvent dans $M \setminus S$, nous avons des champs de Killing qui ont des points fixes
   dans $M \setminus S$ et donc une isotropie non triviale sur $M \setminus S$ : absurde.
   \end{demonstration}
   
   Le champ de vecteurs $X$ obtenu dans la proposition pr\'ec\'edente est de norme constante \'egale ou bien \`a $0$ (si $I$ est un sous-groupe
   \`a un param\`etre unipotent), ou bien \`a $1$ (si $I$ est  sous-groupe \`a un param\`etre  semi-simple).
   
   Nous traitons s\'epar\'ement  ces deux cas  en les deux s\'ections suivantes. 
   
   Avant de clore cette section, remarquons que  le fait (d\'emontr\'e \`a la   proposition~\ref{surface}) que  $S$ soit de dimension $2$ implique que la restriction \`a $S$ de $\cal G$ est un isomorphisme d'alg\`ebres de Lie. En effet, une isom\'etrie locale qui fixe tous les points d'un ouvert de $S$ admet un $1$-jet trivial en chacun de ces points
   et est donc trivale.
   
   \subsection{Isotropie unipotente}
   
   Pla\c cons-nous dans le cas o\`u {\it l'image de $g^{(s)}$ contient au moins une orbite $O_{1}$ de dimension $2$ (contenue dans l'adh\'erence de $O$) dont la composante connexe 
   du stabilisateur est un sous-groupe \`a un param\`etre unipotent de $PSL(2, \RR)$.} 
   
   Consid\'erons un point $u$ de $S$ o\`u le $s$-jet de $g$ appartient \`a l'orbite $O_{1}$.
  Nous allons pr\'eciser la g\'eom\'etrie  de la surface $\Sigma$ passant par $u$ et associ\'ee \`a $O_{1}$  gr\^ace au point iii) de  la  proposition~\ref{surface} :
   
   \begin{proposition}   \label{geodesique}
   
   i) Le champ $X$ est g\'eod\'esique.
   
   ii) La surface $\Sigma$ est totalement g\'eod\'esique et le feuilletage induit par le noyau de $g$ sur $\Sigma$ est transversalement riemannien.
  \end{proposition}
  
  \begin{demonstration}
  
  La m\'etrique lorentzienne $g$ \'etant d\'eg\'en\'er\'ee en restriction \`a $\Sigma$, l'espace tangent \`a $\Sigma$ est le champ de plans $X^{\bot}.$
   Le champ de vecteurs $X$   \'etant tangent \`a $\Sigma$, on peut consid\'erer sa
   d\'eriv\'ee covariante le long de tout champ de vecteurs $W$ contenu dans $X^{\bot}$ (tangent \`a $\Sigma$). Comme la $g$-norme de $X$ est constante, nous avons  d\'ej\`a que pour tout $W \in X^{\bot}$~: $2 \cdot g(\nabla_{W}X,X)= W \cdot g(X,X)=0.$ Le champ de plans $X^{\bot}$ est donc stable par l'op\'erateur $\nabla_{\cdot}X$, que l'on peut interpr\'eter comme une section
   au-dessus de $\Sigma$ du fibr\'e $End(X^{\bot})=End(T \Sigma)$     des endomorphismes de $X^{\bot}$ .  
   
  i)  Constatons d'abord  que le champ $X$ est g\'eod\'esique.  Pour s'assurer que $X$ est bien g\'eod\'esique consid\'erons  l'action du champ de Killing unipotent qui fixe un point $u$ de $\Sigma$ (associ\'e \`a la composante connexe $I$ du stabilisateur de $O_{1}$). Si $H(u) \in X^{\bot}(u)$
   est un vecteur de $g$-norme unitaire, la diff\'erentielle du temps $T$ du flot de champ de Killing envoie le vecteur $H(u)$ sur $H(u) +T \cdot X(u)$. Comme ce flot
   pr\'eserve la connexion $\nabla$ et le champ de vecteurs $X$, il vient que $\nabla_{H(u)}X = \nabla_{H(u) + T \cdot X(u)}X$, ce qui implique que $\nabla_{X}X$ s'annule
   au point $u$.

 Le champ $X$ est alors g\'eod\'esique et,
   par cons\'equent, l'op\'erateur $\nabla_{\cdot}X$ contient le champ $X$ dans son noyau.
   
  ii)  L'autre valeur propre
   de l'op\'erateur $\nabla_{\cdot}X$ (qui est constante sur $\Sigma$ car l'action de $\cal G$ pr\'eserve  $X$ et est  transitive sur $\Sigma$)
   est n\'ecessairement  nulle : dans le cas contraire l'op\'erateur $\nabla_{\cdot}X$ serait diagonalisable et tout \'el\'ement de $\cal G$  devrait pr\'eserver la d\'ecomposition de
   $X^{\bot}$ en deux espaces propres de dimension $1$;  or ceci n'est pas r\'ealis\'e pour notre champ de Killing unipotent dont la diff\'erentielle ne fixe aucune autre droite de $X^{\bot}$ \`a part celle engendr\'ee par $X$. Il reste que le champ 
   d'endomorphismes $\nabla_{\cdot}X$ est nilpotent (d'ordre au plus $2$)~: l'image de $  \nabla_{\cdot}X$ est incluse dans le noyau de   $  \nabla_{\cdot}X$. 
     Deux cas se pr\'esentent : ou bien  $\nabla_{\cdot}X$  est nul,  ou bien    le noyau de $\nabla_{\cdot}X$ et l'image de $\nabla_{\cdot}X$ co\"{\i}ncident avec la
     droite engendr\'ee par $X$ (l'unique droite de $X^{\bot}$ invariante par l'action de $\cal G$). Dans les deux cas le calcul suivant est valide pour tous les
     champs de vecteurs    locaux $W_{1}$ et $W_{2}$ tangents \`a $\Sigma$~:
     $g(\nabla_{W_{1}}   W_{2}, X)= W_{1} \cdot g(W_{2}, X) - g(\nabla_{W_{1}}  X, W_{2})=0,$ le deuxi\`eme terme du membre de droite de l'\'egalit\'e \'etant nul car
    $ \nabla_{W_{1}}  X   $ est contenu dans l'image de $\nabla_{\cdot}X$ et donc colin\'eaire \`a $X$ (tandis que $W_{2} \in X^{\bot}$). Ceci montre que     $  \nabla_{W_{1}}   W_{2}     \in X^{\bot}$, et que
    $\Sigma$ est totalement g\'eod\'esique.
     
     Comme $g$ est d\'eg\'en\'er\'ee en restriction \`a la surface totalement g\'eod\'esique $\Sigma$, le feuilletage engendr\'e sur $\Sigma$ par le champ $g$-isotrope
     $X$ est transversalement riemannien~\cite{Ze}. 
      \end{demonstration}
        
        Analysons  maintenant l'action de l'alg\`ebre de Lie $\cal G$ des champs de Killing de $g$ au voisinage d'un point $u$ de $\Sigma$. La restriction \`a $\Sigma$ de chaque  \'el\'ement de $\cal G$ donne un champ de vecteurs (tangent \`a $\Sigma$) d\'efini dans un  voisinage de $u$ dans $\Sigma$
        dont le flot  pr\'eserve $X$ et donc, en particulier, le feuilletage (transversalement riemannien) $\cal F$ d\'efini par $X$. 
        
        D\'esignons par $\cal H$ l'id\'eal de $\cal G$   form\'e
        par les \'el\'ements de $\cal G$ dont la restriction \`a $\Sigma$ agit trivialement sur la transversale de $\cal F$. 
        Les \'el\'ements de $\cal H$ fixent chaque feuille de  $\cal F$ et ils commutent avec $X$ : ils sont de la forme $f \cdot X$ avec $f$ fonction analytique constante sur les orbites de $X$ (en particulier, les \'el\'ements de $\cal H$ sont $g$-isotropes sur $\Sigma$).  Remarquons que l'alg\`ebre de Lie $\cal H$ est de dimension 
        deux ($\cal H$ ne peut \^etre de dimension $3$ car dans ce cas l'action de $\cal G$ ne serait pas transitive sur $\Sigma$) et le quotient $\cal G$/$\cal H$ qui agit non trivialement sur la transversale de $\cal F$ est n\'ecessairement de dimension $1$ (isomorphe \`a l'alg\`ebre de Lie $\RR$ agissant par translation).  Comme $\cal F$ est transversalement riemannien, le flot de $X$ (comme le flot de tout
        champ de  vecteurs tangent au feuilletage) pr\'eserve la restriction de $g$ \`a $\Sigma$.

        On peut  choisir sur un voisinage ouvert $U$ de $u$ dans $\Sigma$ un champ de vecteurs analytique  $H$ de $g$-norme constante \'egale \`a $1$ et tel que
        $\lbrack X, H \rbrack=0$ (il suffit de d\'efinir $H$ de $g$-norme constante \'egale \`a $1$ sur une petite transversale \`a $\cal F$ et de le transporter par le flot
        du champ  $X$ qui pr\'eserve la restriction de $g$ \`a $\Sigma$). D\'efinissons sur  un voisinage de $u$ dans $\Sigma$ un syst\`eme de  coordonn\'ees $(x,h)$ centr\'e en $u$  et tel que  $\frac{\partial}{\partial x}=X$ et $\frac{\partial}{\partial h}=H$.   Dans ces coordonn\'ees l'expression locale de la  forme quadratique  $g$ restreinte \`a $\Sigma$ est $dh^2$ et
       les \'el\'ements de $\cal H$ restreint \`a $\Sigma$ sont de la forme $f(h)  \frac{\partial}{\partial x}$ (car ils pr\'eservent $\frac{\partial}{\partial x}$
       et $dh^2$).
       Par ailleurs l'alg\`ebre de Lie        $\cal G$/$\cal H$       est engendr\'ee par un champ de vecteurs de la forme $\frac{\partial}{\partial h} +l(h) \frac{\partial}{\partial x},$ 
       pour une certaine fonction analytique $l$ d\'efinie dans un voisinage de l'origine dans $\RR$. Nous avons la relation 
       $\lbrack \frac{\partial }{\partial h} + l(h) \frac{\partial}{\partial x}, f(h) \frac{\partial}{\partial x} \rbrack  = f'(h) \frac{\partial}{\partial x}.$

       Pour comprendre la structure de l'alg\`ebre de Lie $\cal G$  qui agit par isom\'etries affines  pour la restriction de la connexion $\nabla$
       \`a $\Sigma$, nous allons pr\'eciser la structure locale de cette connexion sur $\Sigma$.

       \begin{proposition}
       
     i).   Si $R$ est le tenseur de courbure de $(\Sigma, \nabla )$ alors $R(X,H)X=0$ et $R(X,H)H=\gamma X$, o\`u $\gamma$ est un nombre r\'eel.
     
     ii). La restriction de $\nabla$ \`a $\Sigma$ est localement sym\'etrique. De plus $(\Sigma, \nabla)$ est localement isom\'etrique ou bien \`a la connexion
     canonique  du groupe affine de la droite r\'eelle si $\gamma \neq 0$, ou bien \`a la connexion
     canonique de $\RR^2$ si $\gamma=0$.
     
     \end{proposition}

     Avant de passer \`a la preuve rappelons que  le rev\^etement universel $AG$ du groupe affine de la droite  est un groupe de Lie  de dimension $2$
     qui peut \^etre vu comme l'ensemble des couples $(a,b) \in \RR^2,$ avec $a$ strictement positif,  muni de la multiplication $(a,b) \cdot (a',b')= (aa', ab'+b).$ Ce groupe admet une unique
     connexion lin\'eaire, bi-invariante, sans torsion, compl\`ete et localement sym\'etrique. Le groupe d'isom\'etries de cette connexion est form\'e par
     les translations \`a droite et \`a gauche et il est, par cons\'equent, isomorphe au produit $AG \times AG$~\cite{Zeg}.

     \begin{demonstration}
     Les arguments suivants sont inspir\'es de~\cite{Zeg} (partie 8).
     
     Remarquons que $\nabla_{H}X$ ne d\'epend pas du champ de vecteurs $H$ de $g$-norme unitaire chosi. En effet, si $H'$ est
     un autre champ de vecteurs local tangent \`a $\Sigma$ et de $g$-norme constante \'egale \`a $1$, alors $H'=H +f \cdot X$, pour une certaine
     fonction analytique locale $f$ et comme $X$ est g\'eod\'esique, $\nabla_{H} X= \nabla_{H'}X$. Ce champ de vecteurs est donc invariant
     par l'action de $\cal G$; comme l'action de $\cal G$ est   transitive sur $\Sigma$, ceci  implique qu'il existe un nombre r\'eel $a$ tel que $\nabla_{H}X=aX$.
     
     Comme $\lbrack X, H \rbrack=0$, nous avons que $R(X,H)X= \nabla_{X}  \nabla_{H}  X - \nabla_{H} \nabla_{X}X=\nabla_{X}(aX)-0=0$.
     
     Pour la deuxi\`eme formule sur la courbure remarquons que le terme $R(X,H)H$ ne d\'epend pas du choix de $H$ car si $H'=H +f \cdot X$,
     alors la premi\`ere \'egalit\'e implique que $R(X,H)f \cdot X=0$ et donc $R(X,H)H=R(X,H')H'$. Comme pr\'ec\'edemment,  ceci implique que le champ de vecteurs
     $R(X,H)H$ est invariant par l'action (transitive) de $\cal G$ et, par cons\'equent, $R(X,H)H= \gamma X$, avec $\gamma \in \RR$.  Un calcul direct implique que $\gamma =-a^2$.

      La deuxi\`eme assertion de la proposition est  prouv\'ee dans~\cite{Zeg} sous l'hypoth\`ese suppl\'ementaire $\nabla_{H}H=0$ qui est utilis\'ee
      pour montrer que la courbure est parall\`ele (connexion localement sym\'etrique sur $S$). Nous montrons que
      dans notre situation  nous pouvons nous passer de l'hypoth\`ese faite sur $H$. Remarquons d'abord que $H$ \'etant de $g$-norme constante, $\nabla_{H}H$ est orthogonal \`a
      $H$. Comme $\nabla_{H}H \in X^{\bot}$ ceci implique qu'il existe une fonction analytique $g$ telle que $\nabla_{H}H =g \cdot X$.
      
      Rappelons que la d\'eriv\'ee du tenseur $R$ est donn\'ee par la formule : $$\nabla R(A,B,C,D)=\nabla_{A}(R(B,C)D) -R(\nabla_{A}B,C)D -R(B, \nabla_{A}C)D-R(B,C) \nabla_{A}D,$$ o\`u les   vecteurs $A,B,C, D$ prennent les valeurs $X$ ou $H$. Ce n'est que dans le cas o\`u trois des vecteurs A,B,C,D valent $H$ que l'hypoth\`ese suppl\'ementaire est
      utilis\'ee dans~\cite{Zeg}. Nous traitons ici ce  cas et pour le reste de la preuve nous renvoyons \`a~\cite{Zeg} (proposition 8.4 et proposition 9.2).

1)   $ \nabla R(H,X,H,H)=$ $$ \nabla_{H}R(X,H)H-R(\nabla_{H}X,H)H- R(X, \nabla_{H}H)H- R(X,H) \nabla_{H}X=$$ $$\nabla_{H} \gamma X -R(aX, H)H -R(X,gX)H-R(X,H)aX=
      a \gamma X -a \gamma X =0.$$

2)     $ \nabla R(H,X,H,H)=$ $$ \nabla _{H}R(X,H)H -R(\nabla_{H}X,H)H -R(X, \nabla_{H}H)H -R(X,H)\nabla_{H}H =$$ $$ \nabla_{H} \gamma X - R(aX,H)H - R(X,gX)H - R(X,H)gX=
     a \gamma X- a \gamma X -0 -0 =0.$$
     
 3)    $\nabla R(H,H,H,X)=$ $$ \nabla_{H} R(H,H)X - R( \nabla_{H}H,H)X - R(H, \nabla_{H}H)X-R(H,H) \nabla_{H}X=$$
     $$0-R(gX,H)X -R(H,gX)X-0=0.$$
      \end{demonstration}

       Revenons \`a pr\'esent \`a
      l'alg\`ebre de Lie $\cal G$. Dans le cas o\`u $\gamma=0$ notre alg\`ebre de Lie se plonge dans l'alg\`ebre de Lie du groupe des transformations affines de $\RR^2$ qui
      pr\'eservent la forme quadratique $dh^2$ et un champ de vecteurs isotrope et parall\`ele et elle est engendr\'ee par $\frac{\partial}{\partial x}, \frac{\partial}{\partial h}, h\frac{\partial}{\partial x}$.
        Cette alg\`ebre de Lie est l'alg\`ebre de Lie du groupe Heisenberg car $\frac{\partial}{\partial x}$ est dans le centre et $\frac{\partial}{\partial x}= \lbrack \frac{\partial}{\partial h}, h \frac{\partial}{\partial x} \rbrack$. L'alg\`ebre de Heisenberg  est unimodulaire (car nilpotente) et cette situation a d\'ej\`a \'et\'e analys\'ee et \'elimin\'ee comme contradictoire dans le 
        lemme~\ref{unimodulaire}.      
        
        Supposons maintenant que $\gamma \neq 0$. Dans ce cas $\cal     G$ se plonge dans l'alg\`ebre de Lie du produit $AG \times AG$  (comme
        $AG$ est simplement connexe et complet, toute isom\'etrie locale de $\nabla$ se prolonge en une isom\'etrie globale~\cite{Amo}, \cite{Gro}, \cite{No}). 
        Comme $\cal G$ est de dimension  $3$, un calcul simple montre que, quitte \`a permuter les deux copies du groupe affine,  $\cal G$ est
        n\'ecessairement   engendr\'ee par trois \'el\'ements de l'alg\`ebre de Lie de $AG \times AG$ de la forme $(t,0)$, $(0,t)$ et $(w, \alpha w)$, o\`u $t$, $w$ et $\alpha$ d\'esignent
       respectivement  le g\'en\'erateur infinit\'esimal des translations, le g\'en\'erateur infinit\'esimal des homoth\'eties et un nombre r\'eel. Rappelons que les champs de vecteurs
       $t$ et $w$ sont  invariants par les translations \`a droite et fournissent une base de  l'alg\`ebre de Lie
       du groupe affine telle que  $\lbrack t,w \rbrack =t$. Par d\'efinition, la connexion canonique du groupe affine est d\'efinie par $\nabla_{w}t = \frac{1}{2} \lbrack w,t \rbrack =-\frac{1}{2} t,
       \nabla_{t}  t=   \nabla_{w}  w=0.$ Le tenseur de courbure $R$ de cette connexion est tel que $R(w,t)t=0$ et $R(t,w)w=- \frac{1}{4}t.$

        Traitons d'abord le cas $\alpha=0$. Il vient que $\cal G$ est n\'ecessairement  isomorphe au produit directe de $\RR$ par l'alg\`ebre de Lie du groupe affine et, par cons\'equent,
  $\cal G$ admet un centre non trivial. L'\'el\'ement central est n\'ecessairement de la forme $f(h) \frac{\partial}{\partial x}$ avec $f'(h)=0$,
  ce qui implique que $f$ est constante. Nous avons donc que $X$ est la restriction \`a $S$ d'un champ de Killing local qui se trouve dans le centre de  $\cal G$. Comme ce champ de Killing est invariant par l'action de $G$ sur lui-m\^eme, il fournit un champ de Killing globalement
  d\'efini sur $M \setminus S$ et de norme constante (n\'ecessairement \'egale \`a $0$ car ce champ de Killing se prolonge
  sur un ouvert de $S$ en le champ $g$-isotrope $X$).  La propri\'et\'e de prolongement de champs de Killing implique alors l'existence d'un champ de Killing
  isotrope d\'efini dans  voisinage ouvert $U$ dans $M$ d'un point $u$ de $\Sigma$ et invariant par l'action de $\cal G$. Ce champ est n\'ecessairement  non singulier (voir~\cite{AS}, lemme 6.1 ou~\cite{Zegh}). 
  
  Notre situation est similaire avec celle d\'ecrite dans le lemme~\ref{champ}; ici le contexte est local, mais la m\^eme preuve fonctionne. En effet, la pr\'esence
  du champ de vecteurs isotrope invariant permet de consid\'erer au-dessus de l'ouvert $U$ un sous-fibr\'e principal $\cal R$ du fibr\'e des rep\`eres orthonorm\'ees $R(U)$
  de groupe structural unipotent isomorphe (en tant que groupe alg\'ebrique) \`a  $\left(  \begin{array}{cc}
                                                                 1   &   T \\
                                                                 0     &  1 \\
                                                                 \end{array} \right)$.
                                                                   
    Le $s$-jet $g^{(s)}$ de $g$ s'exprime comme une application
                                                             $\left(  \begin{array}{cc}
                                                                 1   &   T \\
                                                                 0     &  1 \\
                                                                 \end{array} \right)$-\'equivariante de $\cal R$ dans l'espace vectoriel d'une repr\'esentation alg\'ebrique lin\'eaire du groupe unipotent et
                                                                 l'image de cette application est constitu\'ee de l'adh\'erence d'une orbite (qui correspond au $s$-jet de $g$ sur $U \setminus \Sigma$).
         
         Les orbites du groupe unipotent \'etant ferm\'ees, ceci implique que l'image de $g^{(s)}$ est constitu\'ee d'une seule orbite et, par cons\'equent, $g$ est localement
         homog\`ene sur $U$ : absurde.
         
         Il reste \`a \'eliminer le cas $\alpha \neq 0$. Nous utilisons un argument global en consid\'erant la composante connexe de $\Sigma$ dans $S$. Nous d\'esignons \'egalement 
         par $S$ cette composante connexe pour ne pas alourdir les notations. Quitte \`a consid\'erer un rev\^etement fini de $S$, le champ de vecteurs isotrope $X$ construit dans 
         la proposition~\ref{surface}  est d\'efini sur $S$. Comme $S$ est une surface lisse qui admet un champ de vecteurs non singulier, il vient que $S$ est diff\'eomorphe
         \`a un tore. Nous avons d\'emontr\'e \`a la proposition~\ref{geodesique} que $S$ est totalement g\'eod\'esique. La connexion sur $S$ est localement model\'ee sur
         le groupe affine $AG$ (autrement dit, il existe sur $S$ une $(AG \times AG, AG)$-structure) et la restriction de $g$ \`a $S$ s'identifie dans les coordonn\'ees $(a,b)$ du groupe affine (mod\`ele du demi-plan sup\'erieur) \`a une  forme quadratique de rang $1$ proportionnelle \`a $\frac{da^2}{a^2}$. En effet,  cette forme quadratique  est compl\`etement d\'efinie \`a partir de la connexion canonique du groupe
         affine par la formule $R(t,u)u=-<u,u> t$, o\`u $R$ est le tenseur de courbure de la connexion (le noyau de $<,>$ est donn\'e par la direction du champ $t$).
         
        L'expression des champs de vecteurs  $t$ et  $w$ dans le mod\`ele du demi-plan sup\'erieur est respectivement $\frac{\partial}{\partial b}$ et $a\frac{\partial}{\partial a} +
        b \frac{\partial}{\partial b}$.

         Le champ g\'eod\'esique
         isotrope $X$ s'identifie \`a un champ de la forme $k(a) \frac{\partial}{\partial b}$, avec $k$ une fonction analytique,  tandis que les g\'en\'erateurs   $(t,0)$, $(0,t)$ et $(w, \alpha w)$ de l'alg\`ebre de Lie $\cal G$ s'identifient  respectivement
         aux champs de vecteurs  $\frac{\partial}{\partial b}$, $a \frac{\partial}{\partial b}$ et $(\alpha +1) a \frac{\partial}{\partial a} + b \frac{\partial}{\partial b}$.

         Nous avons vu que
         l'action de $\cal G$ pr\'eserve le champ de vecteurs $X$.
         Si tout champ de la forme $k(a) \frac{\partial}{\partial b}$ est invariant par $(t,0)$ et par $(0,t)$, en revanche  l'invariance de $X$ par $(w, \alpha w)$          fournit la condition
         $\lbrack (\alpha +1) a \frac{\partial}{\partial a} + b \frac{\partial}{\partial b}, k(a) \frac{\partial}{\partial b} \rbrack =0$, ce qui donne l'\'equation
         $(\alpha +1 ) a \frac{\partial k}{\partial a}= k$.
         
         Cette \'equation admet des solutions non identiquement nulles si et seulement si $ \alpha \neq  -1$ et dans ce cas la solution g\'en\'erale est $k(a)= Ca^{\frac{1}{\alpha +1}},$ avec $C$ une constante r\'eelle.  Il vient que $\alpha \neq -1$ et que l'expression de $X$ dans le mod\`ele du demi-plan sup\'erieur est de la forme
         $Ca^{\frac{1}{\alpha +1}} \frac {\partial}{\partial b},$ avec $C$ un nombre r\'eel non nul.
         
        D\'emontrons  que la $(AG \times AG, AG)$-structure sur $S$ est compl\`ete.
         D\'esignons par $\tilde S$ le rev\^etement universel de $S$ et par $\tilde X$ l'image r\'eciproque du champ $X$ sur $\tilde S$. L'application d\'eveloppante de notre
          $(AG \times AG, AG)$-structure envoie  $\tilde S$ sur un ouvert de $AG$ et le  champ $\tilde X$ sur le  champ $k(a) \frac{\partial}{\partial b}$.
          Le champ $\tilde X$ \'etant complet (car $X$ est d\'efini sur la vari\'et\'e compacte $S$), pour chaque orbite de $\tilde X$  l'application d\'eveloppante r\'ealise un diff\'eomorphisme
          entre un ouvert connexe de $\tilde S$    contenant l'orbite en question  et invariant par $\tilde X$   et  une bande ouverte  horizontale 
          de la forme $\rbrack a - \epsilon , a + \epsilon \lbrack \times \RR$ dans l'espace des param\`etres $(a,b)$ du groupe affine.         
         
                 Comme $S$ est un tore, on peut construire un champ de vecteurs $H$, globalement d\'efini  sur $S$, et de norme constante \'egale \`a $1$. Pour cela il suffit
          de choisir un champ de vecteurs qui n'est en aucun point colin\'aire \`a $X$ et de le diviser par sa norme (cette norme est non nulle car le noyau de la restriction de $g$ \`a $S$ est engendr\'e par $X$). Comme le feuilletage engendr\'e par $X$ est
          transversalement riemannien et les orbites de $H$ sont parm\`etr\'ees par la logueur, le flot de $H$ envoie n\'ecessairement  une orbite de $X$ sur une autre orbite de $X$ (sans  respecter n\'ecessairement  
          le param\'etrage). L'image r\'eciproque $\tilde H$ de $H$ par le rev\^etement universel de $S$ est un champ de vecteurs complet. Par connexit\'e, une orbite
          du flot de $\tilde H$  intersecte chaque orbite de $\tilde X$. Ceci est suffisant pour conclure \`a la compl\'etude de notre $(AG \times AG, AG)$-structure (voir~\cite{Zeg} proposition 9.3).
          
                   Dans ce cas    le groupe fondamental de $S$ s'identifie \`a  un sous-groupe discret  $\Gamma$ 
         de $AG \times AG$ isomorphe \`a $\ZZ \oplus \ZZ$ qui agit librement, discontinument et proprement sur $AG$. Consid\'erons $p_{1}$ et $p_{2}$ les projections
         de $AG \times AG$ respectivement sur la premi\`ere et sur la seconde coordonn\'ee. Il vient que $p_{1}(\Gamma)$ est un sous-groupe commutatif du groupe
         affine qui est, par cons\'equent, contenu dans un sous-groupe \`a un param\`etre $l_{1}^t$ de $AG$. Dans $AG \times AG$, le sous-groupe \`a un param\`etre
         $(l_{1}^t,1)$ centralise $\Gamma$ et fournit un champ de Killing $K_{1}$ de  $\nabla$ globalement d\'efini sur le quotient $S$.  
         
         En utilisant la projection sur la deuxi\`eme coordonn\'ee de $AG \times AG$, on construit de la m\^eme mani\`ere un deuxi\`eme champ de Killing $K_{2}$
         globalement d\'efini sur $S$. Comme $K_{1}$ et $K_{2}$ sont lin\'eairement ind\'ependants et que $\cal G$ est de codimension $1$ dans l'alg\`ebre de Lie
         des champs de Killing de $(AG, \nabla)$, il existe une combinaison lin\'eaire non triviale de $K_{1}$ et $K_{2}$ qui appartient \`a $\cal G$. On vient donc
         de construire un champ de Killing $K$ de $\cal G$ globalement d\'efini sur $S$. 
         
         Commen\c cons par montrer que $K$ ne peut pas \^etre isotrope. Supposons par l'absurde que $K$ est en tout point isotrope. Alors les champs $K$ et $X$ sont
         globalement d\'efinis sur $S$ et en tout point colin\'aires. La surface $S$ \'etant compacte, le ''quotient'' de $K$ par le champ de vecteurs non singulier $X$ est donc une fonction born\'ee. Il est de m\^eme du quotient des relev\'es $\tilde K$ et $\tilde X$ dans le rev\^etement universel de $\tilde S$. Par ailleurs, $\tilde S$ est identifi\'e \`a $AG$
         et dans ce mod\`ele l'expression de $\tilde X$ est  $Ca^{\frac{1}{\alpha +1}} \frac {\partial}{\partial b},$ tandis que  $\tilde K$ est de la forme $\lambda (t,0) + \mu (0,t)=
         (\mu \cdot a + \lambda ) \frac{\partial}{\partial b},$ avec
         $\lambda$ et $\mu$ des constantes r\'eelles dont au moins une est non nulle.
         Le quotient est donc de la forme $\frac{\mu \cdot a + \lambda }     {Ca^{\frac{1}{\alpha +1}}}$. Comme $\alpha \neq  -1$, cette fonction est born\'ee pour $a  >0$ si et seulement
         si $\lambda=0$ et $\alpha =0$ : ce cas a \'et\'e d\'ej\`a trait\'e.
         
         Supposons maintenant $K$ non isotrope. Dans ce cas on peut supposer $K$  de la forme $(w, \alpha w) + \lambda (t,0) + \mu (0,t),$ avec $\lambda$ et $\mu$ des constantes r\'eelles. Comme $K$ est de norme constante sur $S$, $\nabla_{K}  K$ est orthogonal \`a $K$ et donc en tout point colin\'aire \`a $X$. Comme pr\'ec\'edemment,  le quotient
         des champs de vecteurs, globalement d\'efinis sur $S$, $\nabla_{K}  K$ et $X$ doit \^etre une fonction born\'ee.
         
         Par ailleurs, nous calculons la fonction quotient dans le mod\`ele du demi-plan sup\'erieur. L'expression du champ $K$ dans le mod\`ele est
         $(\alpha +1) a \frac{\partial}{\partial a} + (\mu a + b + \lambda ) \frac{\partial}{\partial b}$. Pour calculer $\nabla_{K}  K$ on d\'ecompose $K$
         dans la base $t$ et $w$ de l'alg\`ebre de Lie (les coefficients de cette d\'ecomposition \'etant des fonctions) et on utilise que les champs $t$ et $w$ sont g\'eod\'esiques, tandis que $\nabla_{w}t =- \frac{1}{2} t. $   
         Il vient que $K= (\alpha +1) w + (\mu a -    \alpha b + \lambda)t$ et un calcul direct donne $\nabla_{K} K=(\mu a -    \alpha b - \lambda \alpha)   t$.
         Le quotient entre ce champ et $X$ est donn\'e par la fonction $\frac{      (\mu a -    \alpha b - \lambda \alpha)  }{Ca^{\frac{1}{\alpha +1}}}$ qui est non  born\'ee
         sur le demi-plan sup\'erieur quelque soit $\alpha \neq 0$  : absurde.
         
              \subsection{Isotropie semi-simple}
   
   Nous consid\'erons ici qu'il {\it existe dans l'image de $g^{(s)}$ une orbite $O_{1}$ de dimension $2$ (contenue dans l'adh\'erence de $O$) et dont la composante connexe du stabilisateur  est un sous-groupe \`a un param\`etre semi-simple de $PSL(2, \RR)$.}

  Comme dans le cas pr\'ec\'edent,  nous consid\'erons un point $u$ de $S$ o\`u le $s$-jet de $g$ appartient \`a $O_{1}$ et la surface $\Sigma$ passant par $u$, associ\'ee \`a l'orbite $O_{1}$ par le point iii) de la proposition~\ref{surface}.

   Nous allons montrer  que dans ce cas il existe un champ de Killing global de $g$-norme constante \'egale \`a $1$. Commen\c cons par la proposition suivante :
   
   \begin{proposition}  \label{centre}
   i) L'alg\`ebre de Lie $\cal G$ admet un centre de dimension $1$.
   
   ii) Il existe un champ de Killing global $K$ sur $M$ invariant par l'action de $\cal G$. En particulier, le pseudo-groupe des isom\'etries locales de $g$ qui pr\'eservent $K$
   agit transitivement sur $M \setminus S$ et $K$ est de $g$-norme constante (\'egale \`a $1$).
   \end{proposition}   
   
   \begin{demonstration}
   
  i) Regardons la restriction de $\cal G$ \`a $\Sigma$ (cette restriction est un isomorphisme d'alg\`ebres de Lie). L'action de $\cal G$ pr\'eserve $\Sigma$, le champ de vecteurs $X$
   tangent \`a $\Sigma$ 
   de $g$-norme constante \'egale \`a $1$ construit dans la section pr\'ec\'edente, ainsi que le feuilletage de dimension un  $\cal F$ donn\'e par le noyau de la restriction de $g$
   \`a $\Sigma$  (c'est aussi le feuilletage obtenu en intersectant l'orthogonal $X^{\bot}$ et le c\^one isotrope de $g$). 
   
   D\'esignons par $\cal H$ la sous-alg\`ebre de Lie de $\cal G$ qui agit trivialement sur la transversale de $\cal F$. Alors $\cal G / \cal H$ agit sur la transversale 
   du feuilletage $\cal F$ en pr\'eservant $g$. Il r\'esulte que $\cal G / \cal H$ est de dimension au plus $1$ et que $\cal H$ est de dimension au moins $2$.  En fait la dimension
   de $\cal H$  est n\'ecessairement \'egale \`a $2$ (si $\cal H$ est suppos\'ee de dimension $3$,  le morphisme d'\'evaluation de $\cal G$ en un point de $\Sigma$ \`a une image de dimension $1$ et l'action de $\cal G$ n'est pas transitive sur $\Sigma$ : absurde.) Les \'el\'ements de $\cal H$ sont tangents au feuilletage $\cal F$ et donc $g$-isotropes
   sur $\Sigma$.
   
    Choisissons
   un point $u \in \Sigma$ : il existe alors un \'el\'ement de $\cal H$ qui ne s'annule pas en $u$ et qui est, par continuit\'e, non nul sur un voisinage $U$ de $u$ dans $\Sigma$. Le flot de ce champ
   pr\'eserve $X$, ce qui implique qu'il existe un syst\`eme de coordonn\'ees centr\'e en $u$ dans lequel le  champ de Killing non singulier
   de $\cal H$ s'exprime  $\frac{\partial}{\partial h}$ et $X = \frac{\partial}{\partial x}$. Dans ces coordonn\'ees, la restriction \`a $\Sigma$ de $g$ est $dx^2$ et, en particulier,
   le feuilletage $\cal F$ est transversalement riemannien. Tout champ de Killing qui s'annule en un point agit trivialement sur la transversale de $\cal F$  et est donc un \'el\'ement de $\cal H$.
 En particulier, un  champ de Killing singulier en $u$ qui engendre le sous-groupe \`a un param\`etre (semi-simple) du groupe d'isotropie en $u$ est n\'ecessairement  un \'el\'ement de $\cal H$ qui pr\'eserve $\frac{\partial}{\partial x}$ et donc de la forme  $f(h) \frac{\partial}{\partial h}$, avec $f$ une fonction analytique s'annulant en $0$. Une base de $\cal H$
 est constitu\'ee des champs $\frac{\partial}{\partial h}$ et $f(h) \frac{\partial}{\partial h}$.

    Le feuilletage $\cal F$ \'etant transversalement riemannien, un champ de Killing dont la valeur au point $u$ est $\frac{\partial}{\partial x}$ est n\'ecessairement  une combinaison
    lin\'eaire de $\frac{\partial}{\partial x}$ et d'un champ de vecteurs tangent au feuilletage (et qui pr\'eserve  X). Un suppl\'ementaire de $\cal H$ dans $\cal G$ est donc engendr\'e par un champ de vecteurs de la forme $\frac{\partial}{\partial x} + l(h) \frac{\partial}{\partial h},$ avec $l$ fonction analytique d\'efinie sur un voisinage de $0$ dans $\RR$. 
    
    Consid\'erons la projection de $\cal G$ sur l'alg\`ebre de Lie des germes de champs de vecteurs d\'efini au voisinage  de l'origine sur l'axe des $h$ (sur une feuille de $\cal F$), parall\`element
    \`a l'axe des $x$ (autrement dit, on ``efface'' la composante selon $\frac{\partial}{\partial x})$.  Cette projection est un morphisme d'alg\`ebres de Lie dont l'image est
    l'alg\`ebre de Lie engendr\'ee par $\cal H$ et  le champ de vecteurs $l(h) \frac{\partial }{\partial h}$.     Cette projection ne peut pas \^etre un isomorphisme :
     sinon $\cal G$ est une alg\`ebre de Lie de dimension $3$ qui agit transitivement  sur un espace de dimension $1$ et d'apr\`es un th\'eor\`eme classique du \`a S. Lie~\cite{Lie}  $\cal G$ est isomorphe \`a l'alg\`ebre de Lie unimodulaire $PSL(2, \RR)$ : absurde.
    
    Il vient que le noyau du morphisme pr\'ec\'edent n'est pas trivial et que $\frac{\partial}{\partial x}$ est un champ de Killing.
    Dans les coordonn\'ees choisies, la restriction \`a $\Sigma$ de $\cal G$ est donc engendr\'ee par les champs $(\frac{\partial}{\partial x}, \frac{\partial}{\partial h}, f(h) \frac{\partial}{\partial h})$. Comme $\frac{\partial}{\partial x}$ est un champ de Killing, nous avons  que $X$ est la restriction \`a $U$ d'un champ de Killing. Comme $X$ commute avec tous les \'el\'ements de $\cal G$,  $X$ repr\'esente  un \'el\'ement  non trivial 
    du centre de  $\cal G$. La restriction de $\cal G$ \`a $\Sigma$ \'etant  un isomorphisme, l'alg\`ebre de Lie engendr\'ee par les champs $(\frac{\partial}{\partial x}, \frac{\partial}{\partial h}, f(h) \frac{\partial}{\partial h})$ est de dimension $3$ : il vient que la fonction $f$ est non constante et que le centre de $\cal G$ est de dimension $1$, engendr\'e par     $\frac{\partial}{\partial x}.$

    ii)   Un \'el\'ement du centre de $\cal G$ fournit un champ de Killing invariant par l'action de $G$ sur lui-m\^eme, ce qui donne un champ de Killing d\'efini sur $M \setminus S$
    et de norme constante. La propri\'et\'e de prolongement de champs de Killing au voisinage d'un point $u$ de $\Sigma$ nous a permis de voir que ce champ de Killing 
    se prolonge sur un voisinage de $u$ dans $M$ et qu'en restriction \`a   $\Sigma$,  ce prolongement  est \'egal  \`a $X$ et donc de norme constante \'egale \`a $1$.
    
    Consid\'erons au voisinage d'un point $m$ de $M \setminus S$ un germe $K$ de champ de Killing qui se trouve dans le centre de $\cal G$ et est de norme constante \'egale \`a $1$.
    La propri\'et\'e de prolongement des champs de Killing et le principe de monodromie permettent de prolonger le germe $K$ le long de tout chemin continu issu de $m$ et contenu
    dans $M$. Le r\'esultat de ce prolongement le long d'un lacet d'origine $m$ et non homotope \`a $0$ est un germe de champs de Killing en $m$ de norme constante \'egale \`a 
    $1$ et qui est \'egalement un \'el\'ement central de $\cal G$ : ce champ co\"{\i}ncide n\'ecessairement  avec le champ initial $K$ (eventuellement sur un rev\^etement double de $M$). On vient de montrer que $K$ se prolonge de mani\`ere univoque sur $M$ en
    un  champ de Killing de norme constante \'egale \`a $1$ et qui est pr\'eserv\'e par l'action de $\cal G$.
    \end{demonstration}

   \begin{lemme}  \label{principal}
   
   Il existe sur $M$ un champ de vecteurs analytique  $\tilde Y$ qui est $g$-isotrope  dont le lieu d'annulation est $S$ et tel que tout champ de Killing local
   de $g$  pr\'eserve n\'ecessairement  $\tilde Y$.
   
   \end{lemme}
   
   \begin{demonstration}
   
   La pr\'esence sur $M$ du champ de vecteurs $K$ de $g$-norme constante \'egale \`a $1$ permet de r\'eduire le groupe structural du fibre des rep\`eres orthonorm\'es  $R(M)$ 
   \`a un sous-groupe semi-simple de la forme  $\left(  \begin{array}{cc}
                                                                 T   &   0\\
                                                                 0     &  T^{-1} \\
                                                                 \end{array} \right),$   avec $T \in \RR^{*}_{+}.$

  Le fibr\'e $R(M)$ poss\`ede alors  un sous-fibr\'e principal $\cal R$
   de groupe structural  $\RR^*$. Le $s$-jet de la m\'etrique lorentzienne donne  en restriction \`a $\cal R$ une application $\RR^*_{+}$-\'equivariante 
  du  fibr\'e principal $\cal R$
   dans  une repr\'esentation de $\RR^*_{+}$ sur un espace vectoriel $\RR^{N}$ de la forme $T \cdot  (x_{1}, \ldots, x_{N})= (T^{2 m_{1} } x_{1}, \ldots T^{2m_{N}}x_{N})$, o\`u les $m_{i}$
   sont des entiers  (on ne consid\`ere que 
   les coordonn\'ees $x_{i}$ sur lesquelles l'image de $\cal R$  est non reduite \`a $0$). Dans la repr\'esentation pr\'ec\'edente l'adh\'erence d'une orbite contient, \`a c\^ot\'e de l'orbite initiale, 
   au plus une autre orbite (qui est n\'ecessairement  de dimension $0$).  L'image de notre application est  alors constitu\'ee d'une orbite
   de dimension $1$ et d'une orbite de dimension $0$ (point fixe de l'action) qui se trouve dans son adh\'erence; l'orbite de dimension $0$ est pr\'ecis\'ement l'image
   de la restriction de $\cal R$ \`a $S$.
   Comme  les orbites ne sont  pas toutes ferm\'ees, il vient que  tous les entiers $m_{i}$ sont de m\^eme signe (ou nuls).    
   
   L'image de $\cal R$ par $g^{(s)}$ \'etant constitu\'ee seulement de deux orbites, il vient que le $s$-jet de $g$ est le m\^eme
   en chaque point de $S$ (il correspond \`a l'orbite de dimension $0$ sous l'action de $\RR^*_{+}$ et \`a l'orbite $O_{1}$ sous l'action de $PSL(2, \RR)$). L'orbite $O_{1}$
   est, par cons\'equent, l'unique orbite, autre que $O$, contenue dans l'image  $g^{(s)}(R(M))$ et le pseudo-groupe des isom\'etries locales agit transitivement sur l'ensemble $S$.
   
    Consid\'erons une des coordonn\'ees $x_{i}$ sur laquelle
   l'action de $\RR^*$ n'est pas triviale (l'entier $m_{i}$ est non nul) et consid\'erons la projection de $(x_{1}, \ldots, x_{N}) \to x_{i}$. Cette application \'etant $\RR^{*}_{+}$-\'equivariante,
   elle fournit une application de $\cal R$ dans $\RR$ qui est \'equivariante pour l'action de $\RR^{*}_{+}$ sur $\RR$ donn\'ee par $T \cdot x_{i}=T^{2m_{i}} x_{i}$ et qui s'annule
   exactement sur la restriction de $\cal R$ \`a $S$. 
   
   Remarquons
   \`a pr\'esent que pour $m_{i}=1$ l'application pr\'ec\'edente s'interpr\`ete comme un champ de vecteurs isotrope $Y$ sur $M$ et qui s'annule sur $S$. En effet : si dans le mod\`ele
   lin\'eaire donn\'e par l'action des $PSL(2, \RR)$ sur les formes quadratiques  de la forme $ax^2 +2bxy +cy^2$, le champ $K$ est repr\'esent\'e par le vecteur $2xy$ (et le groupe 
   $\RR^*_{+}$ par le stabilisateur de $2xy$), alors
   dans la base $(x^2,2xy, -y^2)$ le champ $Y$ s'exprime $(x_{i},0,0)$. Le champ $Y$ est bien isotrope et s'annule exactement sur $S$.

   Pour un $m_{i}$ quelconque,
   notre application s'interpr\`ete comme une section analytique du fibr\'e $TM^{\otimes m_{i}}$ (l'isomorphisme entre $TM$ et $T^{*}M$ fournit par la m\'etrique lorentzienne permet
   de se ramener sans perte de g\'en\'eralit\'e au cas o\`u $m_{i} >0).$ Par construction cette section prend valeurs dans le c\^one des puissances $n$-\`emes des
   vecteurs isotropes et son ensemble d'annulation est $S$.  Comme la multiplication par le r\'eel  positif $T^{2m_{i}}$ pr\'eserve chaque demi-droite de $\RR$, et que l'image de l'application de $\cal R$ dans $\RR$ est l'adh\'erence d'une orbite, il vient que cette image est incluse dans une demi-droite ferm\'ee (de signe constant dans $\RR$). Quitte \`a composer
   la projection pr\'ec\'edente sur l'axe $x_{i}$ avec la fonction $x_{i} \to x_{i}^{\frac{1}{n}}$, si l'image est contenue dans la demi-droite $x_{i} \geq 0$, ou bien avec la fonction $x_{i} \to 
   (-x_{i})^{\frac{1}{n}}$, si
   l'image est contenue dans la demi-droite $x_{i} \leq 0$, on construit un champ de vecteurs $\tilde Y$ tel que $Y = \tilde Y ^{\otimes n}$. Le champ $\tilde Y$ est analytique
   sur $M \setminus S$, mais il n'est pas (encore) clair que ceci reste vrai au voisinage des points de $S$ \`a cause de la non d\'erivabilit\'e de la fonction ''racine $n$-\`eme''
   en $0$.

   Nous d\'emontrons \`a pr\'esent que $\tilde Y$ est bien analytique au voisinage des points de $S$ et, par cons\'equent, dans $M$.

  Remarquons d'abord que $Y$ est invariant par l'action de $\cal G$. En effet, parmi les $1$-jets d'applications
   qui relient  deux points de $M \setminus S$ et qui pr\'eservent  la r\'eduction au groupe structural $\RR^{*}_{+}$ (autrement dit, le champ de vecteurs $K$), ceux qui  pr\'eservent
   aussi le $s$-jet de $g$ (et qui s'int\`egre donc en des isom\'etries locales) sont exactement ceux qui pr\'eservent la section $Y$. Ceci implique que le pseudo-groupe
   des isom\'etries locales de $g$ qui pr\'eservent  $K$ (et qui agit transitivement sur $M \setminus S$) pr\'eserve automatiquement $Y$. L'alg\`ebre des champs de Killing
   de $g$ qui pr\'eservent \'egalement $K$ et $Y$ est donc de dimension $3$. Cette alg\`ebre co\"{\i}ncide alors  avec $\cal G$, ce qui implique que tout champ de Killing
   local de $g$  pr\'eserve n\'ecessairement  $K$ et $Y$.

         Analysons l'action de l'alg\`ebre de Lie des champs de Killing $\cal G$ de $g$ au voisinage d'un point $u$ dans $\Sigma$.  Nous avons vu dans la preuve de la proposition~\ref{centre} qu'il existe un voisinage ouvert $U$ de $u$ dans $\Sigma$ tel que la restriction de $\cal G$ \`a $U$ contient une sous-alg\`ebre $\cal G$$_{1}$ commutative
         et de dimension $2$, engendr\'ee par les champs tangents  $X$ et $H$, qui sont de norme constante \'egale respectivement \`a $1$ et $0$ et non singuliers.
         
          En chaque point de $U$ il existe alors un unique vecteur $Z$ (transverse \`a $U$) uniquement d\'etermin\'e
     par les relations $g(Z)=0, g(X,Z)=0$ et $g(H,Z)=1$ (le champ $Z$ engendre la deuxi\`eme droite isotrope (autre que celle engendr\'ee par $H$) de $X^{\bot}$ et est uniquement d\'etermin\'e sur cette droite par la troisi\`eme relation). Le flot g\'eod\'esique de $Z$ permet de d\'efinir un feuilletage au voisinage de $U$ dont les feuilles sont $exp_{U}(tZ)$,
     \`a $t$ fix\'e. L'action de $\cal G$$_{1}$ sur $U$ pr\'eserve les champs $X$, $H$ et par cons\'equent $Z$. Comme cette  action pr\'eserve $U$ et le champ transverse $Z$, elle  fixe
    chaque feuille $exp_{U}(tZ)$,  \`a $t$ fix\'e.
    
    D\'esignons par $\tilde Z$ l'unique champ g\'eod\'esique (d\'efini dans un voisinage de $u$ dans $M$) qui prolonge $Z$ et par $\tilde X$, $\tilde H$ les champs de vecteurs,
    d\'efinis dans  voisinage de $u$ dans $M$, qui prolongent $X$ et $H$ et qui sont enti\`erement d\'etermin\'es par les relations de commutation $\lbrack \tilde Z, \tilde X \rbrack =
    \lbrack \tilde Z, \tilde H \rbrack =0$. Autrement dit, on transporte $X, H$ par le flot g\'eod\'esique de $Z$. Comme la relation $\lbrack X, H \rbrack =0$ est valable sur
    $U$, il vient que $\lbrack \tilde X, \tilde H \rbrack =0.$ L'action de l'alg\`ebre de Lie $\cal G$$_{1}$ pr\'eserve les champs de vecteurs $\tilde X$, $\tilde H$ et $\tilde Z$.
    
    Consid\'erons un syst\`eme de coordonn\'ees $(x,h,z)$ centr\'e en $u$ et d\'efini au voisinage de $u$ dans $M$ par les relations $\tilde X = \frac{\partial}{\partial x}, 
    \tilde H=\frac{\partial}{\partial h}$ et $\tilde Z = \frac{\partial}{\partial z}.$  Dans ces coordonn\'ees  $\cal G$$_{1}$  est engendr\'ee par les translations
    selon les axes des $x$ et des $h$ (et qui pr\'eservent les plans $z$=constant).

 Rappelons que l'action de $\cal G$ et donc aussi l'action de $\cal G$$_{1}$ pr\'eserve le tenseur $Y$. Il vient que  le tenseur $Y$ est invariant par les translations le long des
 axes $x$ et $y$ et est, par cons\'equent,  enti\`erement  d\'etermin\'e si l'on conna\^{\i}t sa restriction \`a l'axe des $z$; autrement dit, $Y(x,y,z)=Y(0,0,z).$

     D\'esignons par $Y'$ la restriction de $Y$ \`a l'axe des $z$.  Il vient que  $Y'(z)=Y(0,0,z) =f(z)(\bar Y(z))^{\otimes n}$, o\`u $f$ est une fonction analytique  \`a une variable d\'efinie  sur un voisinage de $0$ dans $\RR$ et qui admet un z\'ero isol\'e  en $0$ et $\bar Y(z)$ est un vecteur non nul  et $g$-isotrope de $T_{(0,0,z)} \RR^3$ (cette \'ecriture n'est pas unique).

     Consid\'erons maintenant  $\mu z^l$, avec $\mu$ nombre r\'eel non nul (qu'on peut supposer positif) et $l$ entier strictement positif, le premier jet non nul  de $f$ en $0$. Il vient que le premier jet non nul de $Y'$ \`a l'origine (qui s'interpr\`ete comme un polyn\^ome homog\`ene d\'efini sur la droite $Z(u)$ et \`a valeurs dans l'espace vectoriel $T_{u}M^{\otimes n}$) est de la forme $\mu z^l \bar Y(0)^{\otimes n}= \mu z^l (aX(0)+bH(0) +cZ(0))^{\otimes n}$, avec $a,b,c$ des nombres r\'eels  dont au moins un est non nul. L'invariance de cette expression par le flot du champ de Killing semi-simple $K$ qui fixe $u$
    (flot qui stabilise  la g\'eod\'esique issue de $u$ dans la direction $Z(u)$ sans pr\'eserver le param\'etrage et qui fixe $X(u)$) implique que le premier jet non nul de $f$ \`a l'origine est  de la forme $ \mu z^n$, avec $\mu \in \RR^{*}_{+}$ et $a=b=0$. 
    
    En effet, rappelons que
    l'action de la diff\'erentielle en $u$ du temps $T$ de ce  flot est $T \cdot (X(0), H(0), Z(0))=(X(0),T^2H(0),T^{-2}Z(0))$; l'action sur le $l$-jet de $f$ \`a l'origine sera
    donc $T \cdot \mu z^l=\mu T^{2l}z^l$. 
        Le $l$-jet   de $Y'$ se d\'ecompose en une somme dont les termes sont de la forme
 $\mu z^l \cdot a^pb^qc^rX(0)^p \otimes  H(0) ^q \otimes Z(0)^r,$ o\`u $p,q,r$ sont des entiers positifs dont la  somme vaut $n$. 
 Comme l'action lin\'eaire du temps $T$ du flot de $K$ sur $(T_{0} M)^{\otimes n}$ est diagonalisable dans la base form\'ee par les vecteurs du type 
$X(0)^p \otimes H(0)^q \otimes Z(0)^r,$ et la valeur propre associ\'ee \`a ce vecteur propre est $T^{2(q-r)}$, l'invariance du $l$-jet de $Y'$ par
l'action du flot de $Z$ implique que  tous les termes de la forme $X(0)^p \otimes  H(0)^q \otimes Z(0)^{ r}$ avec
$p+q+r=n$ qui  ont un co\'efficient non 
nul doivent \^etre   associ\'es \`a la m\^eme valeur propre. Comme $\bar Y(0)$ est $g(0)$-isotrope, ceci montre que la seule possibilit\'e est
$a=b=0$ et $l=n$.

    Il vient  qu'au voisinage du point $u$ il existe un champ de vecteurs analytique   $\tilde Y$ tel que $Y=( \tilde Y)^{\otimes n}$ : on peut prendre  $\tilde Y(x,y,z)=l(z) z Y(z),$ o\`u l(z) est la racine $n$-\`eme
    de $\frac{f(z)}{z^n}$ dans un voisinage de $0$ dans $\RR$.

    Cette propri\'et\'e \'etant vraie au voisinage de chaque point de $S$, il r\'esulte que $\tilde Y$ est analytique dans $M$ et que  le tenseur $Y$ est la puissance $n$-\`eme du champ de vecteurs $\tilde Y$. Comme $\cal G$ pr\'eserve $Y$, il vient que $\cal G$ pr\'eserve egalement $\tilde Y$.
  \end{demonstration}

   Il r\'esulte de la preuve pr\'ec\'edente que l'on a :
   
   \begin{proposition}   \label{jet non nul}
   
   Le $1$-jet de $\tilde Y$ en un point de $\Sigma$ est non trivial.
   \end{proposition}
   
   \begin{demonstration}
   
   Comme dans la preuve pr\'ec\'edente consid\'erons un point $u$ de $\Sigma$ et le vecteur $X(u) \in T_{u}\Sigma$ fix\'e par le flot du champ de Killing semi-simple $K$ qui s'annule en
   $u$. Les vecteurs $H(u), Z(u)$ sont les $2$ vecteurs isotropes qui engendrent $X(u)^{\bot}$, $g(H(u), Z(u))=1$  et $H(u) \in T_{u} \Sigma$, tandis que $Z(u)$ est transverse
   \`a $\Sigma$.
   
   D\'esignons par  $Y'$ la restriction de $\tilde Y$ \`a la g\'eod\'esique  issue de $u$ dans la direction $Z(u)$ que l'on param\`etre par $z$. Le champ $Y'$ est analytique et admet un z\'ero (isol\'e) au point $u$ : nous pouvons donc l'exprimer dans un voisinage de $u$ dans la g\'eod\'esique sous la forme $f(z) \bar Y(z)$ et $\bar Y(z)$
   est un champ de vecteurs qui ne s'annule pas sur la g\'eod\'esique. 
   
        Consid\'erons maintenant  $\mu z^l$, avec $\mu$ nombre r\'eel non nul (qu'on peut supposer positif) et $l$ entier strictement positif, le premier jet non nul  de $f$ en $0$. Il vient que le premier jet non nul de $Y'$ \`a l'origine (qui s'interpr\`ete comme un polyn\^ome homog\`ene d\'efini sur la droite $Z(u)$ et \`a valeurs dans l'espace vectoriel $T_{u}M$) est de la forme $\mu z^l \bar Y(0)= \mu z^l (aX(u)+bH(u) +cZ(u))$, avec $a,b,c$ des nombres r\'eels dont au moins un est non nul. L'invariance de cette expression par le flot du champ de Killing semi-simple qui fixe $u$
    (flot qui stabilise  la g\'eod\'esique issue de $u$ dans la direction $Z(u)$ sans pr\'eserver le param\'etrage) implique comme dans la preuve du lemme  pr\'ec\'edent que le premier jet non nul de $f$ \`a l'origine est  de la forme $ \mu z$, avec $\mu$ nombre r\'eel non nul et $a=b=0$. 
    
    Nous avons donc que le $1$- jet de $Y'$ en $u$ vaut  $\mu c z    Z(u)$ et est non nul. En particulier, le $1$- jet de $\tilde Y$ en $u$ est non nul.
     \end{demonstration}

    Nous d\'emontrons maintenant que le champ $\tilde Y$ est n\'ecessairement  g\'eod\'esique.

    \begin{proposition}
    
    Le champ de vecteurs $\tilde Y$ est g\'eod\'esique.
    \end{proposition}

    \begin{demonstration}

    L'action de $\cal G$ pr\'eserve $Y$, la connexion $\nabla$ et donc aussi $\tilde Y$ et $\nabla_{\tilde Y} \tilde Y$. Comme $\cal G$ agit
    transitivement sur $M \setminus S$, le champ $\nabla_{\tilde Y} \tilde Y$ est $g$-norme constante sur $M \setminus S$. Par analyticit\'e, la norme de $\nabla_{\tilde Y} \tilde Y$
    est constante sur $M$.    
        
     Remarquons par ailleurs que $\nabla_{\tilde Y}
\tilde Y$ est un champ de vecteurs  sur $M$ qui s'annule au moins en les points o\`u
$\tilde Y$ s'annule (sur $S$). Il  s'agit donc d'un champ de vecteurs $g$-isotrope. De plus comme $\tilde Y$  
est de $g$-norme constante ($g$-isotrope) on a que  $\nabla_{\tilde Y}\tilde Y$ est orthogonal 
\`a $\tilde Y$. Sur l'ouvert $M \setminus S$ le champ $\tilde Y$ est non singulier et  le champ de plans $\tilde Y^{\bot}$ contient un unique champ de droites 
$g$-isotrope qui est pr\'ecis\'ement engendr\'e par $\tilde Y$.  Les champs de vecteurs  $\nabla_{\tilde Y} \tilde Y$  et
$\tilde Y$ sont alors en tout point colin\'eaires, ce qui implique l'existence d'une fonction m\'eromorphe
$f$ telle que $\nabla_{\tilde Y} \tilde Y=f  \tilde Y$. Comme le pseudo-groupe des isom\'etries locales
de $g$ qui pr\'eservent $Y$ agit transitivement sur l'ouvert $M \setminus S$ et que ce
pseudo-groupe pr\'eserve automatiquement la connexion $\nabla$ et le champ de 
vecteurs $\tilde Y$ on a que la fonction $f$ est constante et il existe donc une constante
r\'eelle  $\lambda$ telle que $\nabla_{\tilde Y} \tilde Y=\lambda  \tilde Y$.
  
  On peut interpr\'eter l'op\'erateur $\nabla_{\cdot} \tilde Y$ comme une section du fibr\'e vectoriel
$T^*M \otimes TM$ des endomorphismes du fibr\'e tangent. Les fonctions sym\'etriques en les
$3$-valeurs propres de la section $\nabla_{\cdot} \tilde Y$  sont des fonctions analytiques et  constantes 
 sur  $M \setminus S$ (action transitive du pseudo-groupe des isom\'etries locales de $g$ qui pr\'eservent $Y$) et elles sont donc constantes sur $M$, ainsi que les valeurs
propres. Une de ces trois valeurs propres est \'egale \`a $\lambda$.

 D\'esignons par 
$div \tilde Y$ la divergence du champ de vecteurs $\tilde Y$ par rapport \`a la forme volume
 $vol$ induite sur $M$ par la m\'etrique lorentzienne : $L_{\tilde Y} vol =div \tilde Y \cdot vol,$ o\`u $L_{\tilde Y}$ est la d\'eriv\'ee de Lie dans la direction de $\tilde Y$. Rappelons
que cette divergence en un point $m$ de $M$ n'est rien d'autre que la trace de l'endomorphisme
$(\nabla_{\cdot} \tilde Y)(m)$.  Comme le pseudo-groupe des isom\'etries locales pr\'eservant $\tilde Y$ agit transitivement sur $M \setminus S$, la fonction analytique divergence de
$\tilde Y$ est constante, \'egale \`a un nombre r\'eel $a$. Notons par $\psi^t$ le temps $t$ du flot de $\tilde Y$ : nous avons alors que $(\psi^t)^*vol=exp(a t) \cdot vol.$ Comme
le flot de $\tilde Y$ doit pr\'eserver $\int_{M}vol$, il vient que $a=0$ et que la divergence de $\tilde Y$ est nulle.

   Nous avons donc que la trace
  de l'op\'erateur $\nabla_{\cdot} \tilde Y$ est nulle. Par ailleurs,  on a vu \`a la proposition~\ref{surface} que $S$ est n\'ecessairement  une surface et comme $\tilde Y$ s'annule sur $S$,   au moins deux valeurs propres de $\nabla_{\cdot} \tilde Y$ sont nulles (le noyau contient l'espace tangent \`a $S$). Il vient  que toutes les trois valeur propres de $\nabla_{\cdot} \tilde Y$ sont nulles.
  Ceci implique que $\lambda =0$ et que $\tilde Y$ est g\'eod\'esique.  
  \end{demonstration}  
  
  Dans la suite nous analysons les jets possibles de $\tilde Y$ en un point d'annulation (situ\'e sur $S$), sachant que ces jets doivent \^etre pr\'eserv\'es par le groupe d'isotropie 
  (semi-simple) de $g$. Le but sera de montrer que ces jets admissibles sont incompatibles avec le caract\`ere g\'eod\'esique de $\tilde Y$. Ceci est r\'ealis\'e dans
  la  proposition  suivante  qui ach\`eve la preuve :

\begin{proposition}   \label{singularite} Consid\'erons  un voisinage ouvert $U$ de l'origine dans $\RR^3$ muni d'une m\'etrique lorentzienne analytique  $g$.

Soit $\tilde Y$ un champ de vecteurs analytique  non identiquement nul $g$-isotrope et qui s'annule \`a l'origine. Supposons qu'il existe un champ (de Killing)
 $K$ dans $U$ qui s'annule \`a l'origine, dont le flot local pr\'eserve \`a la fois $g$ et $\tilde Y$ et dont la diff\'erentielle du flot \`a l'origine est un sous-groupe \`a un param\`etre semi-simple du groupe orthogonal de $(T_{0}U ,g_{0}).$  Alors le $1$-jet du champ $\tilde Y$  \`a l'origine est de la forme $\nu z \frac{\partial}{\partial z}$, avec $\nu \in \RR$ 
\end{proposition}

\begin{demonstration}

Par hypoth\`ese la diff\'erentielle \`a l'origine du flot du champ de Killing $K$  fixe un vecteur $\frac{\partial}{\partial x} \in T_{0}U$ de norme \'egale \`a $1$. Consid\'erons une base $(\frac{\partial}{\partial x},\frac{\partial}{\partial y}, \frac{\partial}{\partial z})$ de l'espace tangent \`a l'origine avec la propri\'et\'e que $\frac{\partial}{\partial x}$ est un vecteur de $g$-norme \'egale \`a $1$, tandis que les vecteurs 
$\frac{\partial}{\partial y}, \frac{\partial}{\partial z}$ sont $g$-isotropes, $g$-orthogonaux \`a $\frac{\partial}{\partial x}$ et $g(\frac{\partial}{\partial y}, \frac{\partial}{\partial z})=1. $
    Dans les coordon\'ees exponentielles fix\'ees par le choix  de la base pr\'ec\'edente le temps $T$ du flot du champ de  Killing $K$ 
     est donn\'e  par la transformation lin\'eaire :
    $T  \cdot (x,y,z)=(x, T^2y, T^{-2}z).$
    
    Un calcul simple donne la forme des champs de vecteurs analytiques  s'annulant \`a l'origine et invariants par ce flot : $Y=f(x,y,z) \frac{\partial}{\partial x}    +
    g(x,y,z) \frac{\partial}{\partial y}    + h(x,y,z) \frac{\partial}{\partial z} $ est invariant par le flot de $K$ si et seulement si $f=F(x,yz), g=yG(x,yz)$ et $h=zH(x,yz)$,
    o\`u $F,G$ et $H$ sont des fonctions analytiques  \`a deux variables s'annulant \`a l'origine.
    
       On en d\'eduit les seuls $1$-jets possibles pour des champs de vecteurs s'annulant \`a l'origine et invariants par le flot de $K$ : un tel $1$-jet 
       est n\'ecessairement  de la forme $Y^{1}= \lambda x \frac{\partial}{\partial x}     + \mu y \frac{\partial}{\partial y}  + \nu z \frac{\partial}{\partial z} $, avec $\lambda, \mu, \nu$,
       des nombres r\'eels eventuellement nuls.
       
       Il s'agit maintenant d'exploiter le fait que $Y$ est $g$-isotrope. Dans les coordonn\'ees exponentielles consid\'er\'ees le $1$-jet \`a l'origine de la m\'etrique lorentzienne
        $g$ vaut $g_{0}=dx^2 +dydz$.  Ceci permet d'en d\'eduire le $2$-jet \`a l'origine de la fonction $g(Y)$~: il s'agit de $g_{0}(Y^{1})=
       2 \mu \nu yz + {\lambda}^2x^2$.
          Comme la fonction  $g(Y)$ est identiquement nulle, ce $2$-jet doit  \^etre     trivial et donc $Y^1$ est  de la forme $\nu z \frac{\partial}{\partial z}$   (ou bien
          de la forme $\mu y \frac{\partial}{\partial y}$ ce qui revient au m\^eme) .
   \end{demonstration}
          
          Remarquons alors que le premier jet non nul de $\nabla_{\tilde Y} \tilde Y$ \`a l'origine 
          vaut $\nabla_{Y^1}  Y^1= \nu^2 z \frac{\partial}{\partial z}$ (en coordonn\'ees exponentielles  la connexion $\nabla$
 est plate \`a l'ordre $0$ et donc ses co\'efficients de Christoffel \`a l'origine sont les m\^emes que ceux de la connexion standard). Comme $\tilde Y$
 est g\'eod\'esique, ce  jet ne peut \^etre nul que si $\nu=0$. Ceci est impossible \`a cause de la proposition~\ref{jet non nul} : absurde.

   {\small

Version abr\'eg\'ee du titre : M\'etriques analytiques lorentziennes de dimension $3$

Mots-cl\'es : vari\'et\'es lorentziennes  - structures rigides- champs de Killing locaux.

Classification math. :  53A55, 53B30, 53C50.

\newpage

{\footnotesize

\vspace{2cm}

Sorin Dumitrescu

\rule{4cm}{.05mm}

D\'epartement de Math\'ematiques d'Orsay 

\'Equipe de Topologie et Dynamique

Bat. 425

\rule{4cm}{.05mm}

U.M.R.   8628  C.N.R.S.

\rule{4cm}{.05mm}

Univ. Paris-Sud (11)

91405 Orsay Cedex

France

\rule{4cm}{.05mm}

Sorin.Dumitrescu@math.u-psud.fr

 \end{document}